\documentclass[a4paper,12pt,reqno]{amsart}
\usepackage{amsmath}
\usepackage{enumerate}
\usepackage{float}
\usepackage{amssymb}
\usepackage{amsfonts}
\usepackage{graphicx}
\usepackage{mathtools} %for prescript
\usepackage[colorlinks]{hyperref}
\renewcommand\eqref[1]{(\ref{#1})} %Need with hyperref
\graphicspath{ {images/} }
\providecommand{\Real}{\mathop{\rm Re}\nolimits}%
\title[Integro-differential equations]{Linear and non-linear integro-differential equations with variable coefficients and Sonin kernels}

\author[S. Nuñez]{Sebastian Nuñez}
\address{
	Sebastian Nuñez:
	\endgraf
	Department of Mathematics
	\endgraf
Center for Research and Advanced Studies 
\endgraf
Av. IPN 2508, 07360, CDMX, Mexico
	\endgraf
	{\it E-mail address} {\rm sebastian.nunez@cinvestav.mx}}

\author[J.E. Restrepo]{Joel E. Restrepo}
\address{
	Joel E. Restrepo:
	\endgraf
	Department of Mathematics
	\endgraf
Center for Research and Advanced Studies 
\endgraf
Av. IPN 2508, 07360, CDMX, Mexico
	\endgraf
	{\it E-mail address} {\rm joel.restrepo@cinvestav.mx}}


\thanks{The paper has been accepted to be published in ``Communications in Nonlinear Science and Numerical Simulation".}

\newtheoremstyle{theorem}%name
{10pt}          % space above
{10pt}  % space below
{\sl}  % bofy font
{\parindent}     % ident - empty=no indent,  \parindent= paragraph indent
{\bf}  % thm head font
{. }    % punctuation after thm head
{ }    % space after thm head: `` ``=normal \newline=linebreak
{}     % thm head specification
\theoremstyle{theorem}

\numberwithin{equation}{section}
\theoremstyle{plain}
\newtheorem{thm}{Theorem}[section]
\newtheorem{prop}[thm]{Proposition}
\newtheorem{cor}[thm]{Corollary}

\theoremstyle{definition}
\newtheorem{defn}[thm]{Definition}
\newtheorem{rem}[thm]{Remark}
\newtheorem{ex}[thm]{Example}

\newtheoremstyle{defi}%name
{10pt}          % space above
{10pt}  % space below
{\rm}  % bofy font
{\parindent}     % ident - empty=no indent,  \parindent= paragraph indent
{\bf}  % thm head font
{. }    % punctuation after thm head
{ }    % space after thm head: `` ``=normal \newline=linebreak
{}     % thm head specification
\theoremstyle{defi}

\begin{document}
 	\begin{abstract}
We study linear and non-linear integro-differential equations, with Sonin and absolutely continuous kernels, and with variable coefficients which are either integrable or continuous functions. In the linear case, for all cases, the unique solution is given explicitly by a convergent infinite series involving compositions of convolutions of the kernels and the variable coefficients. Many examples are provided. For the non-linear case, the existence and uniqueness are established under a suitable condition in the non-linear term. The analytic solution is given by the limit of the iterative sequence found from the method of successive approximations. The exact analytical expressions for several problems are found by computing the iterative sequences that approximate the solution.
\end{abstract}
	\maketitle
%\tableofcontents

\section{Introduction}
The study of the following ordinary differential equations with variable coefficients is a classic topic:
\[
v^{(n)}(t)+a_{n-1}(t)v^{(n-1)}(t)+a_{n-2}(t)v^{(n-2)}(t)+\cdots+a_{1}(t)v^{\prime}(t) = g(t), \,\, t\in[0,T],
\] 
where the variable coefficients $a_i(t)$ $(i=1,\ldots,n)$ can be integrable, continuous, Hölder, of other types of functions. Recent investigations treat similar equations where the local behavior of the derivatives is changed by considering non-local operators. In particular, fractional differential operators can be used. These operators are closely related to the studies of fractional calculus and fractional differential equations (FDEs) \cite{kilbas,oldham,samko}. Applications in these contexts can be found in \cite{hilfer,mainardi}. 

Explicit solutions have been constructed for several classes of linear FDEs with variable coefficients. Approaches varying from  Green's functions \cite{RL,greenfunction}, the Banach fixed point theorem \cite{FRS,analitical}, power series methods \cite{kilbasalpha,vcapl}, Volterra integral equations \cite{vcserbia1,vcserbia2,first}, among others.

In this paper, we follow a classical method that consists in transforming the original integro-differential equation into an equivalent integral equation. Thus, the Banach fixed point theorem and the method of successive approximations can be applied to prove the existence and uniqueness of solutions. For linear FDEs, it seems that the pioneering work in using the latter methodology was \cite{AMS-1938}. This approach has been successfully used in many works, e.g. \cite{kober,transmuted,first,FRS,RL,hilfer-psi,analitical}. Solutions were given by uniformly convergent infinite series involving nested compositions of different types of fractional integrals with multiplier functions.

Linear FDEs with variable coefficients have also been studied by other types of integro-differential operators. For example, fractional calculus with respect to functions \cite{RRS} and operators with Mittag-Leffler kernels \cite{FRS:AB,FRS:P}. 

The studies of FDEs with variable coefficients have opened up the possibility of new works and applications in partial integro-differential equations of wave and heat types \cite{garra, RSade}. Also, some researchers have defined new classes of time-fractional Dirac type operators with time-variable coefficients and with applications in fractional Clifford analysis \cite{faustino,dirac-joel}. 

It is important to highlight that, to the best of our knowledge, the non-linear case of FDEs with variable coefficients is still under review. In this regard, we mention the recent paper of \cite{Li-nonlinear}, where a Green function is constructed to solve it. In our research, we also contributed to the study of the nonlinear case using a different approach. 

Aiming to have a wide class of integro-differential equations with singular kernels, we then consider the following linear (respectively, nonlinear) problem:
\[
(\kappa_{0} * v)'(t) + \sum_{i=1}^{m} \sigma_{i}(t) (\kappa_{i} * v)'(t)  = g(t) \big(\text{or}\quad g(t,v(t))\big), \quad t\in[0,T],     
\]
under the initial condition 
\[
(\kappa_0*v)(0)=a\quad \text{or}\quad v(0)=a,
\]
where $v$ is the unknown function, $\kappa_0$ is a Sonin kernel, $a\in\mathbb{R}$ and $\kappa_i\in AC[0,T]$ for $i=1,\dots,m.$ The variable coefficients and source function are considered in two different function spaces: i) $\sigma_i,g\in L^1[0,T];$ ii) $\sigma_i,g\in C[0,T].$

It is worth pointing out that expressions of the form  $(\kappa_{0} * v)'(t)$ are related to operators of the so-called general fractional calculus \cite{Luchko2021,Yu-sym,Yu-fcaa,Yu-1,Yu-2,Kochubei} and references therein. This idea goes back at least as far as \cite{Kochubei}.

Direct applications of our equation can be found in fractional kinetic models \cite{HauboldMathai2000}, viscoelasticity \cite{Koeller1984,MainardiSpada2011}, and dielectric relaxation \cite{GarrappaMainardiMaione2016}. Sonin kernels play a central role in our approach, as they include classical fractional power types and some other more general.

One of the revelations of this work is that our approach and techniques do not allow us to consider all $\kappa_i's$ as Sonin kernels. Perhaps, the function spaces involved in the solutions have to be improved. We show clearly in our proofs that these kernels must be in $AC[0,T]$ with a regularized condition. Mainly, while finding solutions of the above equation in $L^1[0,T],$ we need to impose more regularity over the kernels to estimate the expressions $(\kappa_i\ast v)'(t).$ This will be discussed in detail in Section \ref{s-integrable}.  
 
The paper is organized as follows: In Section \ref{preli}, we recall all the necessary concepts, elementary results,  function spaces, operators, among others ingredients that will be used throughout the paper. In Sections \ref{s-integrable} and \ref{c-linear}, we study linear integro-differential equations with integrable (respectively, continuous) variable coefficients. Existence, uniqueness, and explicit solutions are provided in spaces $L^{\kappa_0}$ and $C^{\kappa_0}.$ While Sections \ref{nonlinear-integrable} and \ref{nonlinear-continuous} are devoted to the nonlinear case. In the whole paper, several examples are given to show the practicality of our results.

\section{Preliminaries}\label{preli}

\subsection{Basic function spaces}

Let us now introduce the following well-known sets:
\begin{align*}
L^1[0,T]&=\left\{f:[0,T]\to\mathbb{R}\;:\;\big\|f\big\|_{L^1[0,T]}:=\int_0^T\big|f(t)\big|\,\mathrm{d}t<+\infty\right\}; \\
AC[0,T]&=\left\{f:[0,T]\to\mathbb{R}\;:\;f\text{ absolutely continuous on }[0,T]\right\}; \\
AC^n[0,T]&=\left\{f:[0,T]\to\mathbb{R}\;:\;f^{(n-1)}\text{ exists and in }AC[0,T]\right\},\qquad n\in\mathbb{N}.
\end{align*}
Recall that $L^1[0,T]$ is a Banach space according to its norm $\|\cdot\|_{L^1[0,T]}$, and that $AC[0,T]$ is a Banach space with the following norm:
\[
\big\|f\big\|_{AC[0,T]}=|f(0)|+\big\|f'\big\|_{L^{1}[0,T]}.
\]

We continue by proving an elementary result that will be used in our analysis later. 

\begin{prop}\label{C3}
If $f\in L^1[0,T]$ and $v\in C[0,T]$, then $f*v \in C[0,T].$
\end{prop}
\begin{proof}
    Without loss of generality, assume that $t_0<t,$ thus
\begin{align*}
|(f*v)(t)&-(f*v)(t_0)|=\left| \int_0^{t}v(t-s)f(s)\,{\rm d}s - \int_0^{t_0}v(t_0-s)f(s)\,{\rm d}s\right|\\
    &\leqslant \int_0^{t_0}\left(|v(t-s)-v(t_0-s)|\right) |f(s)|\,{\rm d}s+\int_{t_0}^{t}|v(t-s)|\,|f(s)|\,{\rm d}s \\
   &= \int_0^{t_0}|v(t-s)-v(t_0-s)\|f(s)|\,{\rm d}s+\|v\|_{\infty} \int_{t_0}^{t}|f(u)|\, {\rm d}u.
   \end{align*}
Since $f\in L^1[0,T],$ we have $\displaystyle\lim_{t\to t_0}\int_{t_0}^t|f(u)|\, {\rm d}u = 0.$ By Lebesgue's dominated convergence theorem, we obtain 
\[
\lim_{t \to t_0}\int_0^{t_0}|v(t-s)-v(t_0-s)||f(s)|\, {\rm d}s=0.
\]
Hence
\[
\lim_{t\to t_0} |(f*v)(t)-(f*v)(t_0)|= 0,\quad\text{so}\quad f*v \in C[0,T].
\] 
\end{proof}

\subsection{Fractional integro-differential operators}

We start by recalling the fractional integrals and derivatives of Riemann-Liouville types \cite{diethelm,samko}.
\begin{defn}
The Riemann--Liouville fractional integral of order $\beta\in\mathbb{C}$ with $\Real(\beta)>0$ is defined as follows \cite[Sections 2.3 and 2.4]{samko}:
\begin{align*}
^{RL}I_{0}^{\beta}f(t)&=\frac1{\Gamma(\beta)}\int_0^t (t-s)^{\beta-1}f(s)\,\mathrm{d}s,\qquad f\in L^1[0,T].
\end{align*}
The Riemann--Liouville fractional derivative of order $\beta\in\mathbb{C}$ $(\Real(\beta)\geqslant0)$ is given by:
\begin{align*}
\prescript{RL}{}D_0^{\beta}f(t)&=D^{n}\prescript{RL}{}I_0^{n-\beta}f(t) \\ &=\frac1{\Gamma(n-\beta)}\left(\frac{\rm{d}}{\rm{d}t}\right)^{n}\int_0^t (t-s)^{n-\beta-1}f(s)\,\mathrm{d}s,\quad f\in AC^n[0,T],
\end{align*}
where $n:=\lfloor\Real\,\beta\rfloor+1.$
\end{defn}

\subsection{Sonin kernels}

\begin{defn}\cite[Definition 2.1]{Samko1}
A kernel $\kappa \in L^1_{\mathrm{loc}}(\mathbb{R}_+)$ is called a Sonin kernel if there exists a kernel $\tilde{\kappa} \in L^1_{\mathrm{loc}}(\mathbb{R}_+)$ such that
\begin{equation}\label{I0}
(\kappa * \tilde{\kappa})(t)= \int_0^t \kappa(t-s)\tilde{\kappa}(s)\,{\rm d}s=1,\quad\text{for almost all}\quad t \in \mathbb{R}_+.
\end{equation}
The relation \eqref{I0} is symmetric, so $\tilde{\kappa}$ is also a Sonin kernel. As noted in \cite[ p.~2]{Samko1}, both kernels $\kappa$ and $\tilde{\kappa}$ are necessarily unbounded as $t \to 0^{+}$.
\end{defn}

Below are several examples of Sonin kernels. The next example can be found in \cite[p.~5]{Luchko2021}. Consider the following functions:
\begin{equation*}
     h(t)=\frac{t^{-1/2}}{\sqrt{\pi}}, \quad \ell(t)=\frac{t^{-1/2}}{\sqrt{\pi}}, \quad t>0.
\end{equation*}
To verify condition \eqref{I0}, we compute
\begin{align*}
    (h*\ell)(t)&=\int_0^t\frac{(t-s)^{-1/2}s^{-1/2}}{\pi}{\rm d}s=\frac{1}{\pi}\int_0^1(1-u)^{-1/2}u^{-1/2}\,{\rm d}u\\
    &= \frac{1}{\pi}B(1/2,1/2)=\frac{\pi}{\pi}=1,
\end{align*}
where $B(\cdot,\cdot)$ denotes the Beta function. Hence, $(h,\ell)$ forms a Sonin kernel pair.

In general, we can verify whether a pair of functions $(f,g)$ forms a pair of Sonin kernels by using the Laplace transform $\mathcal{L},$ which is defined as $\mathcal{L}\{f(t)\}(s)=F(s)=\int_{0}^{+\infty}e^{-st}f(t)\,{\rm d}t$ for any $\text{Re}(s)>0.$ Indeed, by the convolution property, we know that
\begin{equation*}
\mathcal{L}\{f*g\}(s) = \mathcal{L}\{f\}(s)\, \mathcal{L}\{g\}(s).
\end{equation*}
For $(f, g)$ to be a pair of Sonin kernels, it suffices that their convolution satisfies \eqref{I0} or equivalently 
\[s \, \mathcal{L}\{f\}(s)\, \mathcal{L}\{g\}(s) = 1, \]
provided that their Laplace transforms exist. 

\medskip Let us now present another pair \cite[p.~19]{Samko1}:
\begin{equation*}
f(t) = \frac{\cos(2\sqrt{t})}{\sqrt{\pi t}}, \quad 
g(t) = \frac{\cosh(2\sqrt{t})}{\sqrt{\pi t}}.
\end{equation*}
Their Laplace transforms are
\[
\mathcal{L}\{f\}(s) = \frac{e^{-1/s}}{\sqrt{s}}, \quad 
\mathcal{L}\{g\}(s) = \frac{e^{1/s}}{\sqrt{s}}.
\]
It yields 
\[
s \, \mathcal{L}\{f\}(s)\, \mathcal{L}\{g\}(s) = s \cdot \left(\frac{e^{-1/s}}{\sqrt{s}} \cdot \frac{e^{1/s}}{\sqrt{s}}\right) = 1.
\]
Another interesting example \cite[p.~20]{Samko1} is given by 
\[
\varphi(t) = \frac{e^{-\lambda t}}{\Gamma(\alpha)t^{1-\alpha}}, 
\quad
\phi(t) = \lambda^\alpha +\frac{\alpha}{\Gamma(1-\alpha)} \int_t^{+\infty} \frac{\ e^{-\lambda \tau}}{\tau^{\alpha+1}} {\rm d}\tau,\quad 0\leqslant\alpha\leqslant 1,\quad \lambda\geqslant0.
\]
It is easy to check that
\[
\mathcal{L}\{\varphi\}(s) = \frac{1}{(\lambda+s)^{\alpha}}, \quad 
\mathcal{L}\{\phi\}(s) = \frac{(\lambda+s)^{\alpha}}{s},
\]
and hence
\[
s \, \mathcal{L}\{\varphi\}(s)\, \mathcal{L}\{\phi\}(s) = s\left( \frac{1}{(\lambda+s)^{\alpha}} \cdot \frac{(\lambda+s)^{\alpha}}{s} \right)=1.
\]

\subsection{Special function spaces}

We introduce some spaces that will be used to analyze the solutions of our equations.

\subsubsection{Subspace of $L^1$}

\begin{defn}
Let $\kappa_0$ be a Sonin kernel. We define 
\[
\textit{L}^{\kappa_0}:=\{v \in L^1[0,T] : \kappa_0*v\in AC[0,T]\}
\]
endowed with the norm
\begin{align*}
\|v\|_{L^{\kappa_0}}&=\|v\|_{L^1[0,T]}+ \|\kappa_0*v\|_{AC[0,T]}\\
    &=\|v\|_{L^1[0,T]}+ \|(\kappa_0*v)'\|_{L^1[0,T]}+|(\kappa_0*v)(0)|.
\end{align*}
\end{defn}
Next, we prove that $L^{\kappa_0}$ is a Banach space with the norm $\|\cdot\|_{L^{\kappa_0}}$. 
\begin{thm}\label{BS-i}
The space \(\left(L^{\kappa_0}, \|\cdot\|_{L^{\kappa_0}}\right)\) is a Banach space.
\end{thm}
\begin{proof}
We show that every Cauchy sequence in \(L^{\kappa_0}\) converges to a limit in \(L^{\kappa_0}\). Let \(\{v_n\}_{n=1}^{\infty}\) be a Cauchy sequence in \(L^{\kappa_0}\), then, for every \(\varepsilon > 0\), there exists \(N \in \mathbb{N}\) such that for all \(m, n > N\) we have $\|v_n - v_m\|_{L^{\kappa_0}} < \varepsilon,$ or what is the same
\[
\|v_n - v_m\|_{L^1[0,T]} + \left\|(\kappa_0 * v_n)' -  (\kappa_0 * v_m)'\right\|_{L^1[0,T]} + |(\kappa_0 * v_n)(0) - (\kappa_0 * v_m)(0)| < \varepsilon.
\]
From the above inequality, we can infer that:
\begin{enumerate}[i)]
    \item $\{v_n\}_{n=1}^\infty$ is a Cauchy sequence in \(L^1[0,T]\);
    \item \(\{ (\kappa_0 * v_n)'\}_{n=1}^\infty\) is a Cauchy sequence in \(L^1[0,T]\);
    \item \(\{(\kappa_0 * v_n)(0)\}_{n=1}^\infty\) is a Cauchy sequence in \(\mathbb{R}\).
\end{enumerate}
Since \(L^1[0,T]\) and \(\mathbb{R}\) are Banach spaces, these sequences converge to limits in their respective spaces:
\begin{enumerate}[i')]
    \item There exists \(v \in L^1[0,T]\) such that \(v_n \to v\) in \(L^1[0,T]\);
    \item There exists \(f\in L^1[0,T]\) such that \( (\kappa_0 * v_n)' \to f\) in \(L^1[0,T]\);
    \item There exists \(c \in \mathbb{R}\) such that \((\kappa_0 * v_n)(0) \to c\).
\end{enumerate}
We verify that the candidate limit \(v\) lies in \(L^{\kappa_0}\). We already know \(v \in L^1[0,T]\); it remains to show that \(\kappa_0 * v \in AC[0,T]\). Define \(h_n(t) = (\kappa_0 * v_n)(t)\). Since \(h_n \in AC[0,T]\),
\[
h_n(t) = (\kappa_0 * v_n)(0) + \int_0^t (\kappa_0 * v_n)'(s)\,{\rm d}s.
\]
Define $h(t) := c + \int_0^t f(s)\,{\rm d}s.$ Clearly \(h \in AC[0,T]\), \( h' = f\) a.e. and \(h(0) = c\). We want to prove that \(h_n \to h\) uniformly. In fact, for any \(t \in [0,T]\):
\begin{align*}
|h_n(t) - h(t)| &= \left| (\kappa_0 * v_n)(0) + \int_0^t (\kappa_0 * v_n)'(s)\,{\rm d}s - \left(c + \int_0^t f(s)\,{\rm d}s \right) \right| \\
&\leqslant |(\kappa_0 * v_n)(0) - c| + \int_0^t \left|(\kappa_0 * v_n)'(s) - f(s)\right|\,{\rm d}s \\
&\leqslant |(\kappa_0 * v_n)(0) - c| + \left\| (\kappa_0 * v_n)' - f\right\|_{L^1[0,T]}.
\end{align*}
Both terms on the right-hand side tend to zero as \(n \to +\infty\), so \(\|h_n - h\|_{\infty} \to 0\), i.e., \(h_n \to h\) uniformly. This implies that \(h_n \to h\) in \(L^1[0,T]\). Since \(v_n \to v\) in \(L^1[0,T]\), \(\kappa_0 \in L^1[0,T]\) and Young's convolution inequality, we arrive at
\begin{align*}
\|h_n-\kappa_0 * v\|_{L^1[0,T]}&=\|\kappa_0 * v_n - \kappa_0 * v\|_{L^1[0,T]} \\
&\leqslant \|\kappa_0\|_{L^1[0,T]} \cdot \|v_n - v\|_{L^1[0,T]} \to0,\quad\text{as}\quad n\to+\infty.
\end{align*}
By the uniqueness of the limit in \(L^1\), we conclude that \(h = \kappa_0 * v\) almost everywhere. As \(h \in AC[0,T]\), \(\kappa_0 * v \in AC[0,T]\). Therefore, \(v \in L^{\kappa_0}\). Finally, we just need to prove that \(v_n \to v\) in the \(L^{\kappa_0}\)-norm. Indeed, we have
\[
\|v_n - v\|_{L^{\kappa_0}} = \|v_n - v\|_{L^1[0,T]} + \left\| (\kappa_0 * v_n)' -  (\kappa_0 * v)'\right\|_{L^1[0,T]} +|(\kappa_0 * v_n)(0) - (\kappa_0 * v)(0)|.
\]
We already know \((\kappa_0 * v)'(t) = f(t)\) and \((\kappa_0 * v)(0) = c\). Thus
\[
\|v_n - v\|_{L^{\kappa_0}} = \|v_n - v\|_{L^1[0,T]} + \| (\kappa_0 * v_n)' - f\|_{L^1[0,T]} + |(\kappa_0 * v_n)(0) - c|.
\]
From the above estimates, each term tends to zero as \(n \to +\infty\), therefore \(\|v_n - v\|_{L^{\kappa_0}} \to 0\). We have shown that an arbitrary Cauchy sequence in \(L^{\kappa_0}\) converges to a limit in \(L^{\kappa_0}\) with respect to the \(L^{\kappa_0}\)-norm. Hence, \(L^{\kappa_0}\) is a Banach space.
\end{proof}

\subsubsection{Subspace of $C[0,T]$}

\begin{defn}
Let $\kappa_0$ be a Sonin kernel. Define
\[
C^{\kappa_0}:=\{v \in C[0,T] : \kappa_0*v \in C^1[0,T]\},
\]
endowed with the norm
\begin{align*}
\|v\|_{C^{\kappa_0}}:= &\|v\|_{\infty}+\|\kappa_0*v\|_{C^1[0,T]}\\
=&\|v\|_{\infty}+\|(\kappa_0*v)'\|_{\infty}+\|\kappa_0*v\|_{\infty},
\end{align*}
where $\|\cdot\|_{\infty}$ is the norm in $C[0,T].$ Recall that in $C[0,T],$ the norm $\|\cdot\|_{C[0,T]}$ coincides with $\|\cdot\|_{\infty}.$ 
\end{defn}
Let us prove that $C^{\kappa_0}$ is a Banach space. The following proof is similar to the one given in Theorem \ref{BS-i} and hence we only provide the main steps. 
\begin{thm}
    The space $( C^{\kappa_0},\|\cdot\|_{ C^{\kappa_0}})$ is a Banach space.
\end{thm}
\begin{proof}
Let $\{v_n\}_{n=1}^\infty$ be a Cauchy sequence in $ C^{\kappa_0}$. By definition, we conclude that: $\{v_n\}_{n=1}^\infty$ is a Cauchy sequence in $C[0,T];$ $\{(\kappa_0*v_n)'\}_{n=1}^\infty$ is a Cauchy sequence in $C[0,T];$ $\{\kappa_0*v_n\}_{n=1}^\infty $ is a Cauchy sequence in $C[0,T]$. Since $(C[0,T],\|\cdot\|_{C[0,T]})$ is a Banach space, we have: There exists $v \in C[0,T]$ such that $v_n \to v$ in $C[0,T];$ there exists $g \in C[0,T]$ with $(\kappa_0*v_n)' \to g$ in $C[0,T];$ there exists $h \in C[0,T]$ so that $\kappa_0*v_n \to h$ in $C[0,T].$

Let us see that $v \in C^{\kappa_0}.$ Obviously $v \in C[0,T].$ We now check that $\kappa_0 * v \in C^1[0,T].$ First, observe that $\kappa_0 * v_n \to \kappa_0 * v$ in $C[0,T].$ From Young's convolution inequality,  
\[
\|\kappa_0 * v_n - \kappa_0 * v\|_{\infty}
    \leqslant \|v_n - v\|_{\infty} \|\kappa_0\|_{L^1[0,T]} \to 0, \quad \text{as} \quad n \to +\infty.
\]
By the uniqueness of the limit, $h = \kappa_0 * v$. From $\kappa_0*v_n \in C^1[0,T],$ it follows 
\begin{align*}
    (\kappa_0*v_n)(t)&=(\kappa_0*v_n)(0)+\int_0^t(\kappa_0*v_n)'(s)\,{\rm d}s.
    \end{align*}
Each $v_n$ is continuous, so $(\kappa_0 * v_n)(0) = 0.$ Moreover, from the estimate 
\[
\left|(\kappa_0* v_n)'(t)\right|\leqslant \left\|(\kappa_0*v_n)'\right\|_{\infty}\leqslant M<+\infty,
\]
and by Lebesgue's dominated convergence theorem, we obtain
\begin{equation*}
   (\kappa_0*v)(t)= \lim_{n\to+\infty}(\kappa_0*v_n)(t)=\lim_{n\to+ \infty}\int_0^t(\kappa_0*v_n)'(s)\,{\rm d}s=\int_0^tg(s) \,{\rm d}s.
\end{equation*}
Therefore, $\kappa_0 * v \in C^1[0,T]$ and $(\kappa_0 * v)'(t) = g(t)$ for every $t \in [0,T]$. Consequently, $v \in C^{\kappa_0}.$ Finally, we prove that $v_n \to v$ in the $C^{\kappa_0}$-norm. Indeed, having in mind that $h = \kappa_0 * v$ and $(\kappa_0 * v)' = g,$ we have
\begin{align*}
    \|v_n-v\|_{ C^{\kappa_0}}&=\|v_n-v\|_{\infty}+\|(\kappa_0*v_n)'-(\kappa_0*v)'\|_{\infty}+\|\kappa_0*v_n-\kappa_0*v\|_{\infty} \\
    &=\|v_n-v\|_{\infty}+\|(\kappa_0*v_n)'-g\|_{\infty}+\|\kappa_0*v_n-h\|_{\infty}. 
\end{align*}
Each term tends to zero as $n \to +\infty$, hence $\|v_n - v\|_{ C^{\kappa_0}} \to 0$. 
\end{proof}

\section{Linear integro-differential equations with integrable variable coefficients}\label{s-integrable}

In this section, we study the following linear integro-differential equation with integrable variable coefficients:  
\begin{equation}\label{DE}
(\kappa_{0} * v)'(t) + \sum_{i=1}^{m} \sigma_{i}(t)(\kappa_{i} * v)'(t)  = g(t), \quad t\in[0,T],     \end{equation}
under the initial condition \begin{equation}\label{hic}
(\kappa_0*v)(0)=a,
\end{equation}
where $v$ is the unknown function, $a\in \mathbb{R},$ $\kappa_0$ is a Sonin kernel, $\sigma_i,g\in L^1[0,T]$ and $\kappa_i\in AC[0,T]$ such that $\kappa_i(0)=0$ for $i=1,\dots,m.$   

The main objective in what follows is to prove the existence, uniqueness, and explicit solution of problem \eqref{DE}-\eqref{hic} in the space $L^{\kappa_0}$. We will make the following assumption:
\begin{equation}\label{A1}
{C:=\sup_{s\in[0,T]}\sum_{i=1}^{m}\int_{s}^T|\sigma_i(t)|\,|(\kappa'_i\ast\tilde\kappa_0)(t-s)|\, {\rm d}t< 1.}
 \end{equation} 
 
\medskip The main result of this section is as follows.

\begin{thm}\label{thm-integrable}
    If $\sigma_i,g\in L^1[0,T], \ \kappa_i\in AC[0,T]$ with $\kappa_i(0)=0$ for $i=1,\dots,m$, and assumption \eqref{A1} holds, then problem \eqref{DE}-\eqref{hic} has a unique solution $v\in L^{\kappa_0}$, and is given by the following infinite series convergent in $L^{\kappa_0}$: 
    $$v(t)=\sum_{k=0}^{+\infty}(-1)^k \Lambda^k\big((\tilde{\kappa}_0*g)+a\tilde{\kappa}_0)\big)(t),$$
    with 
    \[
    \Lambda f:=\tilde{\kappa}_0 * \left( \displaystyle\sum_{i=1}^m \sigma_i(t)\, (\kappa_i' * f) \right),
    \]
    where $\Lambda^0 := \mathrm{Id}$ and $\Lambda^k$ $(k \geqslant 2)$ is the $k$-fold composition of $\Lambda$.
\end{thm}
\begin{proof}
The procedure is to find a solution in $L^{\kappa_0}$. First, we transform equation \eqref{DE} with condition \eqref{hic} into an equivalent integral equation, more simple to solve. Hence, for $v \in L^{\kappa_0}$ being a solution of \eqref{DE}-\eqref{hic}, define $u(t):= (\kappa_0 * v)'(t),$ that evidently belongs to $L^1[0,T].$ By \eqref{hic}, 
\[
Iu(t):=\int_0^t u(s)\,{\rm d}s= \int_{0}^{t} (\kappa_0 * v)'(s)\,{\rm d}s = (\kappa_0 * v)(t)-(\kappa_0*v)(0)=(\kappa_0 * v)(t)-a.
\]
By the Sonin condition $(\tilde{\kappa}_0\ast \kappa_0)(t)=1,$ it follows that $(\tilde{\kappa}_0*Iu)(t)+aI\tilde\kappa_0(t)=Iv(t)$ and $v(t)= (\tilde{\kappa}_0 * Iu)'(t)+a\tilde{\kappa}_0(t).$ By \cite[Corollary 7.4]{gri}, we obtain
\begin{align*}
(\tilde{\kappa}_0*Iu)'(t)=Iu(0)\tilde{\kappa}_0(t)+\int_0^t\frac{{\rm d}}{{\rm d}t} Iu(t-s)\tilde{\kappa}_0(s)\,{\rm d}s =(u*\tilde{\kappa}_0)(t),
\end{align*}
i.e. $v(t)=(\tilde{\kappa}_0*u)(t)+a\tilde{\kappa}_0(t).$ Also 
$(\kappa_i * v)'(t) = \big(\kappa_i * [(\tilde{\kappa}_0 * u)+a\tilde{\kappa}_0] \big)'(t).$ According to the above equalities, equation \eqref{DE} can be rewritten as 
\begin{equation}\label{IE}
u(t) + \sum_{i=1}^{m} \sigma_{i}(t) \big(\kappa_i*[(\tilde{\kappa}_0 * u)+a\tilde{\kappa}_0]\big)'(t) = g(t).
\end{equation} 
We have proved: If $v\in L^{\kappa_0}$ is a solution of \eqref{DE}-\eqref{hic} then $u\in L^1[0,T]$ is a solution of \eqref{IE}.

On the other hand, let $u\in L^1[0,T]$ be a solution of \eqref{IE}. Therefore,
\[
(\tilde{\kappa}_0 * u)(t)+a\tilde{\kappa}_0(t) + \tilde{\kappa}_0 *  \sum_{i=1}^m \sigma_i(t)\big(\kappa_i*[(\tilde{\kappa}_0 * u)+a\tilde{\kappa}_0]\big)'(t)=( \tilde{\kappa}_0 * g)(t)+a\tilde{\kappa}_0(t).
\]
Set $v(t)=(\tilde{\kappa}_0 * u)(t)+a\tilde{\kappa}_0(t)$. It yields that
\[
v(t) + \tilde {\kappa} _0 *  \sum_{i=1}^m \sigma_i (t)\big(\kappa_i *v\big)'(t) =  (\tilde {\kappa} _0 * g)(t)+a\tilde{\kappa}_0(t),
\]
which implies
\[
( \kappa _0 * v)(t) +  \kappa _0 * \tilde {\kappa }_0 *  \sum_{i=1}^m \sigma_i (t) (\kappa_i * v)'(t) =  \big(\kappa _0 *\tilde {\kappa}_0 * g\big)(t)+a(\kappa_0*\tilde{\kappa}_0)(t). 
\]
By Young's convolution inequality, $v\in L^1[0,T].$ Moreover $(\kappa_0*v)(t)=a+Iu(t).$ Hence $u(t)=(\kappa_0 * v)'(t),$ $u\in L^1[0,T]$ and $(\kappa_0*v)(0)=a,$ i.e \eqref{hic} holds. For this reason $\kappa_0 * v \in AC[0,T].$ From $(\kappa_0*\tilde{\kappa}_0)(t)=1,$ one has
\[
(\kappa _0 * v)'(t) + \left( 1 * \sum_{i=1}^m \sigma_i (t) (\kappa_i * v)'(t)\right)' =(1 * g)'(t),
\]
which is equivalent to equation \eqref{DE}.    

The task is now to prove the existence and uniqueness of the solution of \eqref{IE}. Define the operator $T:L^1[0,T]\to L^1[0,T]$ as:
\[
Tu(t):=g(t)-\sum_{i=1}^{m} \sigma_{i}(t)\big(\kappa_i*[(\tilde{\kappa}_0 * u)+a\tilde{\kappa}_0]\big)'(t).
\]
\eqref{IE} makes it obvious that $Tu=u.$ Using \cite[Corollary 7.4]{gri} and $\kappa_i(0)=0$ $(i=1,\ldots,m)$, we find that 
\begin{align*}
\big(\kappa_i*[(\tilde{\kappa}_0*u)+a\tilde{\kappa}_0]\big)'(t)&=\kappa_i(0)[(\tilde{\kappa}_0*u)+a\tilde{\kappa}_0](t)+\big(\kappa'_i\ast[(\tilde{\kappa}_0*u)+a\tilde{\kappa}_0]\big)(t)\\
   & =(\kappa'_i*[(\tilde{\kappa}_0*u)+a\tilde{\kappa}_0])(t)=((\kappa_{i,0}*u)+a\kappa_{i,0})(t), 
\end{align*}
where  $\kappa_{i,0}:=\kappa'_i*\tilde\kappa_{0}.$ Hence, we can rewrite the operator $T$ as follows:
\[
Tu(t)=g(t)-\sum_{i=1}^{m} \sigma_{i}(t)\big((\kappa_{i,0} *u)+a\kappa_{i,0}\big)(t).
\]
Next, we see that $T$ is a contraction mapping. In fact, for $u_1,u_2 \in L^1[0,T]$ and any fixed $t \in (0,T),$ we have
\[
Tu_1(t)-Tu_2(t)=-\sum_{i=1}^m\sigma_i(t)\big(\kappa_{i,0} *(u_1-u_2)\big)(t).
\]
By \eqref{A1}, note that
\begin{align*}
  \|Tu_1 - Tu_2\|_{L^1[0,T]}&=\int_0^T\left|\sum_{i=1}^m\sigma_i(t)(\kappa_{i,0}*(u_1-u_2))(t)\right|\,{\rm d}t\\
    &\leqslant \sum_{i=1}^m\int_0^T|\sigma_i(t) |\int_0^t|\kappa_{i,0}(t-s)|\, |u_1(s)-u_2(s)|\,{\rm d}s \, {\rm d}t\\
&=\sum_{i=1}^m\int_0^T|u_1(s)-u_2(s)|\left(\int_{s}^T|\sigma_i(t)|\,|\kappa_{i,0}(t-s)|\, {\rm d}t\right){\rm d}s \\
   & \leqslant C\|\,u_1-u_2\|_{L^1[0,T]}.
    \end{align*}
The mapping $T$ is a contraction. From the Banach fixed-point theorem, there exists a unique point $u\in L^1[0,T]$ such that $Tu=u$. Furthermore, we can find $u$ with the following approximation sequence:
    \begin{align*}
        u_0(t)&=g(t), \\
        u_n(t)&= g(t) - \sum _{i=1}^m\sigma_i(t)\big((\kappa_{i,0} *u_{n-1})+a\kappa_{i,0}\big)(t), \quad n\geqslant1,
    \end{align*}
which converges with respect to $\|\cdot\|_{L^1[0,T]}$ to a limit $u$ in $L^1[0,T].$ Also, by the Young's convolution inequality and $\tilde{\kappa}_0\in L^1[0,T]$, it is straightforward to see that 
\begin{equation}\label{000}
\|\tilde{\kappa}_0*u_n-\tilde{\kappa}_0*u\|_{L^1[0,T]}\to0,\quad\text{as}\quad n\to+\infty.    
\end{equation}
We write this new sequence as follows:
    \begin{align*}
        (\tilde{\kappa}_0*u_0)(t)+a\tilde{\kappa}_0(t) &=(\tilde{\kappa}_0*g)(t)+a\tilde{\kappa}_0(t),\\
        (\tilde{\kappa}_0*u_n)(t)+a\tilde{\kappa}_0(t)&=(\tilde{\kappa}_0*g)(t)+a\tilde{\kappa}_0(t)-\tilde{\kappa}_0*\sum _{i=1}^m\sigma_i(t)\big(\kappa'_{i}*[(\tilde{\kappa}_0*u_{n-1})+a\tilde{\kappa}_0]\big)(t).
    \end{align*}
Denoting $v_n(t)=(\tilde \kappa_0*u_n)(t)+a\tilde{\kappa}_0(t)$ for all $n$, we can rewrite the above sequence as
    \begin{align}
        v_0(t) &= (\tilde{\kappa}_0*g)(t)+a\tilde{\kappa}_0(t), \nonumber\\
        v_n(t)&= v_0(t)- \tilde \kappa_0*\sum _{i=1}^m\sigma_i(t)(\kappa'_i*v_{n-1})(t),\quad n\geqslant1.\label{rel-v}
    \end{align} 
What is left is to show that  $v_0 \in L^{\kappa_0}$ and $v_n\in L^{\kappa_0}$ for any $n\geqslant1.$ Recall that $u_n\in L^1[0,T]$ for any $n\geqslant1$ and $\tilde{\kappa}_0\in L^1[0,T].$ Young's convolution inequality gives
\[
\|v_n\|_{L^1[0,T]}
= \|\tilde{\kappa}_0 * u_n+a\tilde{\kappa}_0\|_{L^1[0,T]}
\leqslant \|\tilde{\kappa}_0\|_{L^1[0,T]}\,\left( \|u_n\|_{L^1[0,T]}+|a|\right) <+\infty,
\]
which means $v_n\in L^1[0,T].$ Of course, $v_0\in L^1[0,T].$ Since $(\kappa_0*v_0)(t)=a+Ig(t),$ we see that $\kappa_0*v_0 \in AC[0,T]$, i.e. $v_0 \in L^{\kappa_0}$. We know that $(\kappa_0*v_n)(t)=a+Iu_n(t),$ so $\kappa_0 * v_n\in AC[0,T]$ and $v_n \in L^{\kappa_0}.$ Observe that $v(t)=(\tilde{\kappa}_0*u)(t)+a\tilde{\kappa}_0(t)$ is in $L^{\kappa_0}.$ Let us check that $v_n\to v$ in $L^{\kappa_0}.$ Since $(\kappa_0*v_n)(t)=a+Iu_n$ and $(\kappa_0*v)(t)=a+Iu,$ we compute
    \begin{align*}
        \| v_n-v\|_{L^{\kappa_0}} &=\|v_n-v\|_{L^1[0,T]}+\|\kappa_0*(v_n-v)\|_{AC[0,T]}\\
        &=\|v_n-v\|_{L^1[0,T]}+\|I(u_n-u)\|_{AC[0,T]}\\
        &=\|v_n-v\|_{L^1[0,T]}+\|u_n-u\|_{L^1[0,T]}\to0,\quad\text{as}\quad n\to+\infty,
    \end{align*} 
where we have used \eqref{000} and $u_n\to u$ in $L^1[0,T].$ As a consequence, the limit function $v$ is the unique solution in $L^{\kappa_0}$ of \eqref{DE}-\eqref{hic}.   

In the following, we construct the explicit solution $v$ using the recurrence relation \eqref{rel-v}. The first approximation is given by
\begin{align*}
    v_1(t)&=v_0(t) - \tilde{\kappa}_0*\sum_{i=1}^m\sigma_i(t)(\kappa_i'*v_0)(t)=\sum_{k=0}^1(-1)^k\Lambda^k(v_0)(t).
\end{align*}
The second approximation reads as 
\begin{align*}
    v_2(t)&=v_0(t)-\tilde{\kappa}_0*\sum_{i=1}^m\sigma_i(t)(\kappa'_i*v_1)(t)\\
          &=v_0(t)-\tilde{\kappa}_0*\sum_{i=1}^m\sigma_i(t)\left(\kappa'_i*\left[v_0(t) - \tilde{\kappa}_0*\left(\sum_{i=1}^m\sigma_i(t)(\kappa_i'*v_0)(t)\right)\right]\right)\\
          &=v_0(t)-\tilde{\kappa}_0*\sum_{i=1}^m\sigma_i(t)(\kappa'_i*v_0)(t)+\tilde{\kappa}_0*\sum_{i=1}^m\sigma_i(t)\left[\kappa'_i*\tilde{\kappa}_0*\sum_{i=1}^m\sigma_i(t)(\kappa_i'*v_0)(t)\right]\\
&=\sum_{k=0}^2(-1)^k\Lambda^k(v_0)(t).
\end{align*}
By induction, we find the $n$-th approximation for all $n\in \mathbb{N},$ that is expressed by \begin{equation}\label{Suc1}
v_n(t)=\sum_{k=0}^n(-1)^k \Lambda^k(v_0)(t)=\sum_{k=0}^n(-1)^k \Lambda^k\big((\tilde{\kappa}_0*g)+a\tilde{\kappa}_0\big)(t).
\end{equation}
Finally, the solution $v$ can be presented as $$v(t)=\lim_{n\to+\infty}v_n(t)=\sum_{k=0}^{+\infty}(-1)^k \Lambda^k\big((\tilde{\kappa}_0*g)+a\tilde{\kappa}_0)\big)(t),$$ where the limit is taken in the norm $\|\cdot\|_{L^{\kappa_0}}$ and the infinite series converges there.
\end{proof}

\begin{rem}
The estimate \eqref{A1} is the most suitable in the current approach. If we want to find another explicit bound, then we end up imposing additional conditions on the kernels $\kappa_i.$ It can be seen from the proof of Theorem \ref{thm-integrable} that
\begin{align*}
  \|Tu_1 - Tu_2\|_{L^1[0,T]}&=\int_0^T\left|\sum_{i=1}^m\sigma_i(t)(\kappa_{i,0}*(u_1-u_2))(t)\right|\,{\rm d}t\\
    &\leqslant \sum_{i=1}^m\int_0^T|\sigma_i(t) |\int_0^t|\kappa_{i,0}(t-s)|\, |u_1(s)-u_2(s)|\,{\rm d}s \, {\rm d}t\\
&=\sum_{i=1}^m\int_0^T|u_1(s)-u_2(s)|\int_{s}^T|\sigma_i(t)|\,|\kappa_{i,0}(t-s)|\, {\rm d}t \, {\rm d}s \\
   & \leqslant \sum_{i=1}^m |\kappa'_{i}(T)|\int_0^T|u_1(s)-u_2(s)|\,{\rm d}s\int_{0}^T|\sigma_i(t)|\,{\rm d}t\int_0^{T} |\tilde{\kappa}_0(\tau)|\, {\rm d}\tau\\
    &=\left(\sum_{i=1}^m |\kappa'_{i}(T)|\, \|\sigma_i\|_{L^1[0,T]} \,\|\tilde{\kappa}_0\|_{L^1[0,T]}\right)\|\,u_1-u_2\|_{L^1[0,T]},
    \end{align*}    
where we need to assume that $|\kappa_i^\prime|$ is increasing. Therefore, we could have a more explicit condition with an extra assumption as follows:
\[
\sum_{i=1}^m|\kappa'_{i}(T)| \, \|\sigma_i\|_{L^1[0,T]} \,\|\tilde{\kappa}_0\|_{L^1[0,T]} < 1.
\]
\end{rem}

\subsection*{Constant coefficients case.} If $\sigma_i(t)=b_i\in\mathbb{R}$ for $i=1,\ldots,m$ in problem \eqref{DE}--\eqref{hic}, the same fixed-point argument applies without assuming $\kappa_i(0)=0$. Let us announce the result without giving a proof, since it follows the same steps as that of Theorem~\ref{thm-integrable}. The only difference is in the contraction estimate, since we now allow $\kappa_i(0)\neq 0$.

\begin{cor}\label{cor:constant-coefficients}
If $\sigma_i(t)=b_i\in\mathbb{R}$, $g\in L^1[0,T]$ and $\kappa_i\in AC[0,T]$ for $i=1,\ldots,m$. Assume that 
\begin{equation*}\label{equ:Asum-Int}
\|\tilde{\kappa}_0\|_{L^1[0,T]}\sum_{i=1}^{m}|b_i|
\left(|\kappa_i(0)|+\|\kappa_i'\|_{L^1[0,T]}\right)<1.
\end{equation*}
Then problem \eqref{DE}--\eqref{hic} has a unique solution $v\in L^{\kappa_0}$, given by 
\[
v(t)=\sum_{k=0}^{+\infty}(-1)^k \varLambda^k(\tilde{\kappa}_0*g+a\tilde{\kappa}_0
)(t),
\]
where the series converges in $L^{\kappa_0}$ and
\[
\varLambda f:=\tilde{\kappa}_0* \left(\sum_{i=1}^{m}b_i(\kappa_i*f)'\right),
\]
 with $\varLambda^0 := \mathrm{Id}$ and $\varLambda^k$ $(k \geqslant 2)$ is the $k$-fold composition of $\varLambda$.
\end{cor}

In the following, we provide an example where the convolution of two integrable functions is not zero, as $t\to0^+.$  
\begin{ex}\label{Example}
Suppose that 
\begin{equation*}
f(t)=t^{\alpha-1}, 
\quad 
\kappa_0(t)=\frac{t^{-\alpha}}{\Gamma(1-\alpha)}, \quad t\in(0,T),\quad 0<\alpha<1.
\end{equation*}
Then
\begin{equation*}
(f*\kappa_0)(t)
=\frac{1}{\Gamma(1-\alpha)}\int_{0}^{1}(1-u)^{\alpha-1}u^{-\alpha}\,{\rm d}u
=\frac{1}{\Gamma(1-\alpha)}\,B(1-\alpha,\alpha).
\end{equation*}
Using the identity $B(x,y)=\dfrac{\Gamma(x)\Gamma(y)}{\Gamma(x+y)}$, we obtain
\begin{equation*}
(f*\kappa_0)(t)
=\frac{1}{\Gamma(1-\alpha)}\cdot 
\frac{\Gamma(1-\alpha)\Gamma(\alpha)}{\Gamma(1)}
=\Gamma(\alpha)
\end{equation*}
and $\displaystyle\lim_{t\to 0^+}(f*\kappa_0)(t)=\Gamma(\alpha)\neq 0.$
\end{ex}

\begin{rem} 
In Theorem \ref{thm-integrable}, a natural question is: Can we guarantee the existence of a solution $v$ in $L^{\kappa_0}$ of problem \eqref{DE} under the initial condition $v(0)=b\in\mathbb{R}$? In this case, our model-spaces, methodology and strategy do not allow to consider it. Let us see that this is not possible. Suppose to the contrary that there exists a solution $v$ in $L^{\kappa_0}$ of \eqref{DE} with $v(0)=b.$ Define $u(t)=(\kappa_0 * v)'(t).$ It is clear that $u\in L^1[0,T]$ and
\[
Iu(t)= \int_{0}^{t} (\kappa_0 * v)'(s)\,{\rm d}s= (\kappa_0 * v)(t)-(\kappa_0*v)(0).
\]
By Example \ref{Example}, in general, $(\kappa_0 * v)(0) \neq 0.$ Let $(\kappa_0 * v)(0) = c$ with $c\in\mathbb{R}$. Since $(\tilde{\kappa}_0 \ast \kappa_0)(t) = 1$, we compute
\begin{align*}
    (\tilde{\kappa}_0*Iu)(t)&=\big(\tilde{\kappa}_0 * ({\kappa}_0 * v)\big)(t)-(\tilde{\kappa}_0*c)(t)= Iv(t)-cI\tilde{\kappa}_0,
\end{align*} 
and thus $v(t)= (\tilde{\kappa}_0 * Iu)'(t)+c\tilde{\kappa}_0$. By \cite[Corollary 7.4]{gri}, one gets
\begin{align*}
\big(\tilde{\kappa}_0 * Iu\big)'(t)=Iu(0)\tilde{\kappa}_0(t)+\int_0^t\frac{{\rm d}}{{\rm d}t} Iu(t-s)\tilde{\kappa}_0(s)\, {\rm d}s=(u*\tilde{\kappa}_0)(t),
\end{align*}
i.e., $v(t) = (\tilde{\kappa}_0 * u)(t) + c\,\tilde{\kappa}_0(t)$. 
Since $\tilde{\kappa}_0$ is a Sonin kernel, it is unbounded as $t\to0^+$. This implies that
$$v(0)=(\tilde{\kappa}_0 * u)(0) + c\,\tilde{\kappa}_0(0)=\pm\infty,$$ 
which contradicts the initial condition $v(0)=b$.
\end{rem}

We present an illustrative example of Theorem \ref{thm-integrable}. The explicit solution is found by estimating the infinite series of the solution.

\begin{ex}
Consider the Sonin kernel $\kappa_0(t)=t^{-\alpha}/\Gamma(1-\alpha)$ for $\alpha\in(0,1).$ Set $\sigma_1(t)=1/{2\beta}$ for $\beta>1$, $\kappa_1(t)=t^\beta/\Gamma(\beta)$ and $g(t)=1.$ Equation \eqref{DE} becomes
\begin{equation}\label{ex1}
\left( \ \frac{t^{-\alpha}}{\Gamma(1-\alpha)}  * v \right)'(t)+\frac{1}{2\beta }\left(\frac{t^\beta}{\Gamma(\beta)}*v\right)'(t)=1, \quad t\in[0,1],
\end{equation}
under the initial condition $\left( \frac{t^{-\alpha}}{\Gamma(1-\alpha)} * v \right)(0)=0.$ We know that $\tilde{\kappa}_0(t)=t^{\alpha-1}/\Gamma(\alpha)$ and
\begin{align*}
\sum_{i=1}^1\int_{s}^1\frac{1}{2\beta}\,&\left|\left(\frac{\beta t^{\beta-1}}{\Gamma(\beta)}  * \frac{t^{\alpha-1}}{\Gamma(\alpha)}\right)(t-s)\right|\, {\rm d}t  \leqslant \frac{1}{2}\int_s^1\frac{(t-s)^{\beta+\alpha-1}}{\Gamma(\beta+\alpha)}\,{\rm d}t \\
&
= \frac{1}{2}\frac{(1-s)^{\beta+\alpha}}{\Gamma(\beta+\alpha)(\beta+\alpha)}\leqslant\frac{1}{2\Gamma(\beta+\alpha+1)} < 1,
\end{align*}
this means that condition \eqref{A1} is satisfied. By Theorem \ref{thm-integrable}, equation \eqref{ex1} has a unique solution $v\in L^{t^{-\alpha}/\Gamma(1-\alpha)}$ and it is represented as
\[
v(t) = \sum_{k=0}^{+\infty} (-1)^k \Lambda^k \left(\frac{t^{\alpha-1}}{\Gamma(\alpha)}*1\right)(t)=\sum_{k=0}^{+\infty} (-1)^k \Lambda^k \left(\frac{t^{\alpha}}{\Gamma(\alpha+1)}\right).
\]
By the semigroup property of the Riemann-Liouville fractional integral, we have
\begin{equation*}
\Lambda f = \tilde{\kappa}_0 * \left( \sum_{i=1}^{1} \sigma_i [\kappa_i' * f] \right) = \frac{t^{\alpha-1}}{\Gamma(\alpha)} * \left( \frac{ t^{\beta-1}}{2 \, \Gamma(\beta)} * f \right) = \frac{1}{2} \, {}^{RL}I_{0}^{\alpha} \, {}^{RL}I_{0}^{\beta}f = \frac{1}{2} \, {}^{RL}I_{0}^{\alpha+\beta}f.
\end{equation*}
Consequently, again by the semigroup property, the $k$-th iteration of the operator satisfies
\begin{equation*}
\Lambda^k f = \left(\frac{1}{2}\right)^k{}^{RL}I_{0}^{k(\alpha+\beta)}f.
\end{equation*}
Substituting the above expressions into the series, the solution turns to
\begin{align*}
v(t) &= \sum_{k=0}^{+\infty} \left(-\frac{1}{2}\right)^k \, {}^{RL}I_{0}^{k(\alpha+\beta)} \left( \frac{t^\alpha}{\Gamma(1+\alpha)} \right) = \sum_{k=0}^{+\infty} \left(-\frac{1}{2}\right)^k
\frac{t^{\alpha+k(\alpha+\beta)}}{\Gamma(\alpha + k(\alpha+\beta) + 1)} \\
&= t^\alpha \sum_{k=0}^{+\infty} \frac{(-\frac{1}{2}t^{\alpha+\beta})^k}{\Gamma(\alpha + k(\alpha+\beta) + 1)}= t^\alpha E_{\alpha+\beta, \alpha+1}\left(-\frac{1}{2}t^{\alpha+\beta}\right).
\end{align*}
So $v(0)=0.$ By \cite[Formula 4.4.5]{Gorenflo}, we obtain
\[
{}^{RL}I_{0}^{1-\alpha}v(t)={}^{RL}I_{0}^{1-\alpha}\left(t^\alpha E_{\beta+\alpha,1+\alpha}\left(-\frac{1}{2} t^{\alpha+\beta}\right)\right)=t\,E_{\beta+\alpha,2}\left(-\frac{1}{2} t^{\beta+\alpha}\right),
\]
which implies that ${}^{RL}I_{0}^{1-\alpha}v(0)=0.$ 

Alternatively, we can rewrite \eqref{ex1} as follows
\begin{equation}\label{ex1.1}
   ^{RL}D^\alpha_0v(t)+ \frac{1}{2}{}^{RL}I^\beta_0 v(t)=1.
\end{equation}
Let us find the solution of \eqref{ex1.1} directly. Applying the Laplace transform to \eqref{ex1.1} gives
\[
s^{\alpha}\mathcal{L}\{v\}(s)-{}^{RL}I_{0}^{1-\alpha}v(0)+\frac{1}{2}s^{-\beta}\mathcal{L}\{v\}(s)=\frac{1}{s},\quad s>0.
\]
Hence
\begin{align*}
\mathcal{L}\{v\}(s)=\frac{1}{s}\frac{1}{s^\alpha+\frac{1}{2s^\beta}}=\frac{1}{s^{1+\alpha}}\frac{1}{1+\frac{1}{2s^{\beta+\alpha}}}.
\end{align*}
By \cite[Theorem 2.3]{Restrepo}, it follows that 
\[
{v}(t)=t^\alpha E_{\beta+\alpha,1+\alpha}\left(-\frac{1}{2} t^{\alpha+\beta}\right),
\]
which coincides with the solution constructed through iteration. 
\end{ex}

\begin{ex}[A fractional Cole--Cole model]
The Cole--Cole dielectric relaxation model, studied in \cite[Equation 3.9]{GarrappaMainardiMaione2016}, is given by
\begin{equation}\label{eq:cole_cole}
{}^{RL}D_0^\alpha\phi_{CC}(t) =-\frac{1}{\tau_*^\alpha}\phi_{CC}(t), \qquad 0<t<T,
\end{equation}
where $0<\alpha<1$, $\tau_*>0$ is the reference relaxation time and $\phi_{CC}(t)$ is the dielectric response function, subject to the initial condition
\[
\left({}^{RL}I_0^{1-\alpha}\phi_{CC}\right)(0)=\frac{1}{\tau_*^\alpha}.
\]
Hence, equation \eqref{eq:cole_cole} can be rewritten as
\[
\left(\frac{t^{-\alpha}}{\Gamma(1-\alpha)}*\phi_{CC}\right)'(t)+\frac{1}{\tau_*^\alpha}(1*\phi_{CC})'(t)=0,
\]
with the initial condition
\[
\left(\frac{t^{-\alpha}}{\Gamma(1-\alpha)}*\phi_{CC}\right)(0)=\frac{1}{\tau_*^\alpha}.
\]
Thus, $\kappa_0(t)=t^{-\alpha}/\Gamma(1-\alpha)$, $\kappa_1(t)=1$, $b_1=1/\tau_*^\alpha$, $g(t)=0$. Since the associated Sonin kernel is $\tilde{\kappa}_0(t)=t^{\alpha-1}/\Gamma(\alpha)$, and given that $\kappa_1(0)=1$ and $\kappa_1'(t)=0$, we obtain
\[
K=\|\tilde{\kappa}_0\|_{L^1[0,T]}|b_1|\left(|\kappa_1(0)|+\|\kappa_1'\|_{L^1[0,T]}\right)=\frac{T^\alpha} {\tau_*^\alpha\Gamma(\alpha+1)}.
\]
Under the assumption $T^\alpha<\tau_*^\alpha\Gamma(\alpha+1),$ i.e. $K<1$, Corollary \ref{cor:constant-coefficients} establishes the existence of a unique solution $\phi_{CC} \in L^{\kappa_0}.$ The operator $\varLambda$ can be expressed as
\[\varLambda w=\frac{t^{\alpha-1}}{\Gamma(\alpha)}*\left(\frac{1}{\tau_*^\alpha}(1*w)'\right)=\frac{1}{\tau_*^\alpha}\,{}^{RL}I_0^\alpha w.
\]
By the semigroup property of the Riemann--Liouville fractional integrals,
\[
\varLambda^k w=\frac{1}{\tau_*^{k\alpha}}{}^{RL}I_0^{k\alpha}w.
\]
Therefore, the solution is
\begin{align*}
\phi_{CC}(t)&=\sum_{k=0}^{+\infty}(-1)^k\varLambda^k(a\tilde{\kappa}_0)(t)=\sum_{k=0}^{+\infty}\left(-\frac{1}{\tau_*^\alpha}\right)^k
{}^{RL}I_0^{k\alpha}\left(\frac{t^{\alpha-1}}{\tau_*^\alpha\Gamma(\alpha)}\right)\\
&=\frac{t^{\alpha-1}}{\tau_*^\alpha}\sum_{k=0}^{+\infty}\frac{(-t^\alpha/\tau_*^\alpha)^k}{\Gamma((k+1)\alpha)}=\frac{t^{\alpha-1}}{\tau_*^\alpha}E_{\alpha,\alpha}\left(-\frac{t^\alpha}{\tau_*^\alpha}\right).
\end{align*}
Let us verify that the initial condition is satisfied:
\begin{align*}
\left(\frac{t^{-\alpha}}{\Gamma(1-\alpha)}*\phi_{CC}\right)(t)&={}^{RL}I_0^{1-\alpha}\phi_{CC}(t)=\frac{1}{\tau_*^\alpha} \sum_{k=0}^{+\infty} \left(-\frac{1}{\tau_*^\alpha}\right)^k{}^{RL}I_0^{k\alpha+1-\alpha}\left(\frac{t^{\alpha-1}}{\Gamma(\alpha)}\right)\\
&=\frac{1}{\tau_*^\alpha}\sum_{k=0}^{+\infty}\frac{(-t^\alpha/\tau_*^\alpha)^k}{\Gamma(k\alpha+1)}=\frac{1}{\tau_*^\alpha}E_\alpha\left(-\frac{t^\alpha}{\tau_*^\alpha}\right).
\end{align*}
Since $E_\alpha(0)=1$, we conclude that $\left({}^{RL}I_0^{1-\alpha}\phi_{CC}\right)(0)
={1}/{\tau_*^\alpha}.$
\end{ex}

\section{Nonlinear integro-differential equations with integrable variable coefficients}\label{nonlinear-integrable}

The purpose of this section is to study a nonlinear version of the equation \eqref{DE}. We consider
\begin{equation}\label{NI1}
 (\kappa_0 * v)'(t)
+ \sum_{i=1}^{m} \sigma_i(t)\,(\kappa_i * v)'(t)
= g(t,v(t)), \quad t \in [0,T],
\end{equation} 
with the initial condition 
\begin{equation}\label{NI2}
    (\kappa_0 * v)(0)=a, \quad a\in \mathbb{R}.
\end{equation}
We assume that $\kappa_0$ is a Sonin kernel, $\sigma_i \in L^1[0,T]$ and $\kappa_i \in AC[0,T]$ with $\kappa_i(0)=0$ for $i=1,\dots,m,$ . Let $\Delta\subseteq\mathbb{R}$ be an open set and let $g:[0,T]\times\Delta\to\mathbb{R}$ be an integrable function satisfying
\begin{equation}\label{NI3}
|g(t,v_1(t))-g(t,v_2(t))| \leqslant \lambda|v_1(t)-v_2(t)|, \quad  t\in[0,T], \quad v_1(t),v_2(t)\in\Delta,
\end{equation}
for some constant $\lambda>0.$ 

In order to obtain nonlinear results, it is necessary to put the following restriction:
\begin{equation}\label{NI4}
D:=\sup_{s\in[0,T]}\left\{\lambda\int_s^T|\tilde{\kappa}_0(t-s)| \, {\rm d}t+\sum_{i=1}^m\int_s^T|\sigma_i(t)||(\kappa'_i*\tilde\kappa_{0})(t-s)|\,{\rm d}t\right\}< 1,
\end{equation}
where $\lambda$ is the constant of \eqref{NI3}.

 \begin{thm}\label{non-linear-thm-int}
       Suppose that $\sigma_i \in L^1[0,T]$, $\kappa_i \in AC[0,T]$ with $\kappa_i(0)=0$ for $i=1,\dots,m,$ and $g$ is integrable on $[0,T]\times\Delta$. If conditions \eqref{NI2}-\eqref{NI3}-\eqref{NI4} are satisfied, then equation \eqref{NI1} admits a unique solution $v \in   L^{\kappa_0}$.
 \end{thm}
 
\begin{proof}
Using similar arguments as in the first part of the proof of Theorem \ref{thm-integrable}, a solution $v \in L^{\kappa_0}$ of \eqref{NI1}-\eqref{NI2} can be equivalent transformed by the change $u(t)=(\kappa_0 * v)'(t)\in L^1[0,T]$ to be a solution of the following integral equation
\begin{equation}\label{NI5}
u(t) + \sum_{i=1}^{m} \sigma_{i}(t) (\kappa_i*[(\tilde{\kappa}_0*u)+a\tilde{\kappa}_0])'(t) = g\big(t,[(\tilde{\kappa}_0*u)+a\tilde{\kappa}_0](t)\big).
\end{equation}
Let us check that \eqref{NI5} admits a solution. Let 
$\widetilde{T}:L^1[0,T]\to L^1[0,T]$ be such that
\begin{equation*}
  \widetilde{T}u(t):=  g\big(t,[(\tilde{\kappa}_0*u)+a\tilde{\kappa}_0](t)\big)-\sum_{i=1}^m\sigma_i(t)(\kappa_i*[(\tilde{\kappa}_0*u)+a\tilde{\kappa}_0])'(t).
\end{equation*}
\eqref{NI5} asserts that $\widetilde{T}u(t)=u(t).$ By \cite[Corollary 7.4]{gri} and $\kappa_i(0)=0$ $(i=1,\ldots,m)$, it yields 
\begin{align*}
   \big(\kappa_i*[(\tilde{\kappa}_0*u)+a\tilde{\kappa}_0]\big)'(t)=(\kappa'_i*[(\tilde{\kappa}_0*u)+a\tilde{\kappa}_0])(t)=((\kappa_{i,0}*u)+a\kappa_{i,0}))(t), 
\end{align*}
where  $\kappa_{i,0}=\kappa'_i*\tilde\kappa_{0}.$ As a result, the operator $\widetilde{T}$ can be expressed equivalently as:
\[ \widetilde{T}u(t)=  g\big(t,[(\tilde{\kappa}_0*u)+a\tilde{\kappa}_0](t)\big)-\sum_{i=1}^m\sigma_i(t)((\kappa_{i,0}*u)+a\kappa_{i,0})(t).\]
We will show that $\widetilde{T}$ is a contraction on $L^1[0,T].$ For $u_1, u_2 \in L^1[0,T]$ we have
\begin{align*}
\widetilde{T}u_1(t)-\widetilde{T}u_2(t)&=g(t,[(\tilde{\kappa}_0*u_1)+a\tilde{\kappa}_0](t))-g(t,[(\tilde{\kappa}_0*u_2)+a\tilde{\kappa}_0](t)) \\
&\hspace{2cm}-\sum_{i=1}^m\sigma_i(t)(\kappa_{i,0}*(u_1-u_2))(t).
\end{align*}
 By the above estimates and \eqref{NI3}, it follows that 
\begin{align*}
 \|\widetilde{T}u_1 &- \widetilde{T}u_2\|_{L^1[0,T]} \leqslant   \int_0^T\left|g(t,[(\tilde{\kappa}_0*u_1)+a\tilde{\kappa}_0](t))-g(t,[(\tilde{\kappa}_0*u_2)+a\tilde{\kappa}_0](t))\right|{\rm d}t\\
&\hspace{2cm}+\int_0^T\left|\sum_{i=1}^m\sigma_i(t)\left[(\kappa_{i,0}*(u_1-u_2))(t)\right]\right| {\rm d}t\\
 &\leqslant \lambda\int_0^T|(\tilde{\kappa}_0*u_1)(t)-(\tilde{\kappa}_0*u_2)(t)| \, {\rm d}t +\sum_{i=1}^m\int_0^T|\sigma_i(t)||(\kappa_{i,0}*(u_1-u_2))(t)| \, {\rm d}t\\
  &\leqslant \lambda\int_0^T \int_0^t|\tilde{\kappa}_0(t-s)||(u_1-u_2)(s)|\,{\rm d}s \,{\rm d}t \\
  &\hspace{2cm}+\sum_{i=1}^m \int_0^T|\sigma_i(t)| \int_0^t|\kappa_{i,0}(t-s)|\,|(u_1-u_2)(s)|\,{\rm d}s\,{\rm d}t\\
 &= \lambda \int_0^T |(u_1-u_2)(s)|\int_s^T|\tilde{\kappa}_0(t-s)|\,{\rm d}t\,{\rm d}s\\
 &\hspace{2cm}+\sum_{i=1}^m\int_0^T|(u_1-u_2)(s)|\int_s^T|\sigma_i(t)||\kappa_{i,0}(t-s)|\,{\rm d}t\,{\rm d}s\\ 
 &=\int_0^T|u_1-u_2|\left(\lambda\int_s^T|\tilde{\kappa}_0(t-s)| \, {\rm d}t+\sum_{i=1}^m\int_s^T|\sigma_i(t)||\kappa_{i,0}(t-s)|\,{\rm d}t\right)\,{\rm d}s\\
&\leqslant D\|u_1-u_2\|_{L^1[0,T]},
 \end{align*}
where $D$ is given from \eqref{NI4}. Finally, $\widetilde{T}$ is a contraction and, by the Banach fixed-point theorem, there exists a unique function $u\in L^1[0,T]$ with $\widetilde{T}u=u$. Consequently, equation \eqref{NI5} admits a unique solution. From the equivalence established between equations \eqref{NI1} and \eqref{NI5}, we conclude that the original problem \eqref{NI1} with conditions \eqref{NI2}-\eqref{NI3}-\eqref{NI4} admits a unique solution $v \in L^{\kappa_0}$.
\end{proof}

\begin{rem}\label{int-nonlinar-rem}[Approximation of the solution] From Theorem \ref{non-linear-thm-int}, we can find $u$ with the following approximation sequence
    \begin{align*}
        u_0(t)&=g(t,0), \\
        u_n(t)&= g\big(t,[(\tilde{\kappa}_0*u_{n-1})+a\tilde{\kappa}_0](t)\big) - \sum _{i=1}^m\sigma_i(t)(\kappa_{i,0}*u_{n-1}+a\kappa_{i,0})(t), \quad n\geqslant1,
    \end{align*}
which converges in the norm $\|\cdot\|_{L^1[0,T]}$ to some function 
$u\in L^1[0,T]$. For $v_n(t)=(\tilde{\kappa}_0*u_n)(t)+a\tilde{\kappa}_0(t),$ the previous sequence takes the following form
    \begin{align}
        v_0(t) &= (\tilde{\kappa}_0*g(\cdot,0))(t)+a\tilde{\kappa}_0(t), \nonumber\\
        v_n(t)&= (\tilde{\kappa}_0*g(\cdot,v_{n-1}))(t)+a\tilde{\kappa}_0(t)- \tilde \kappa_0*\sum _{i=1}^m\sigma_i(t)(\kappa'_i*v_{n-1})(t),\quad n\geqslant1.\label{NI7}
    \end{align} 
It can be easily seen that $\{v_n\}_{n\geqslant0}$ converges to $v=\tilde{\kappa}_0*u+a\tilde{\kappa}_0$ in $L^{\kappa_0}$. So, $v$ is the unique solution in $L^{\kappa_0}$ of \eqref{NI1}.
\end{rem}
We show an application of Theorem \ref{non-linear-thm-int}. In particular, we use Remark \ref{int-nonlinar-rem} to construct the solution iteratively.
 \begin{ex}
 Let $\kappa_0(t)=t^{-\alpha}/\Gamma(1-\alpha)$ and $g(t,v(t))= v(t)/4+t^{-\alpha}$ for $\alpha\in(0,1).$ Consider $\sigma_1(t)=1/8$ and $\kappa_1(t)=t.$  Then equation \eqref{NI1} takes the form
\begin{equation}\label{exam1}
 \,\left( \ \frac{t^{-\alpha}}{\Gamma(1-\alpha)}  * v \right)'(t)+\frac{1}{8}\left(t*v\right)'(t)= \frac{1}{4}v(t)+t^{-\alpha}, \quad t \in [0,1],
\end{equation}
with $\left(\frac{t^{-\alpha}}{\Gamma(1-\alpha)}*v\right)(0)=0.$ We know that
\[\left|\frac{v_1(t)}{4}+t^{-\alpha}-\frac{v_2(t)}{4}-t^{-\alpha}\right|=\frac{1}{4}|v_1(t)-v_2(t)|,\]
i.e. \eqref{NI3} holds. As $\tilde{\kappa}_0(t)= \frac{t^{\alpha-1}}{\Gamma(\alpha)},$ we have
\begin{align*}
          \frac{1}{4}\int_s^1\frac{(t-s)^{\alpha-1}}{\Gamma(\alpha)}{\rm d}t+\sum_{i=1}^1\int_s^1\frac{1}{8}\frac{(t-s)^{\alpha}}{\Gamma(\alpha+1)}{\rm d}t\leqslant\frac{1}{4\Gamma(\alpha+1)}+\frac{1}{8\Gamma(\alpha+2)}<1,
\end{align*}
which implies \eqref{NI4}. By Theorem \ref{non-linear-thm-int}, equation \eqref{exam1} has a unique solution $v \in L^{\frac{t^{-\alpha}}{\Gamma(1-\alpha)}}$. Substituting the above functions into \eqref{NI7}, we obtain the following iterative sequence
\begin{align*}
    v_0(t)&=\left(\frac{t^{\alpha-1}}{\Gamma(\alpha)}*t^{-\alpha}\right)(t),\\
    v_{n}(t)&=\left(\frac{t^{\alpha-1}}{\Gamma(\alpha)}*\left( \frac{v_{n-1}}{4}+t^{-\alpha}\right)\right)(t)-\left(\frac{t^{\alpha-1}}{\Gamma(\alpha)}*\frac{1}{8}
    (1*v_{n-1})\right)(t).
\end{align*}
The previous sequence can be rewritten as
\begin{align*}
    v_0(t)&={}^{RL}I^{\alpha}_0 t^{-\alpha},\\
    v_n(t)&=\frac{1}{4}{}^{RL}I^{\alpha}_0 v_{n-1}(t)+{}^{RL}I^{\alpha}_0 t^{-\alpha}-\frac{1}{8}{}^{RL}I^{\alpha+1}_0 v_{n-1}(t).
\end{align*}
If we iterate the sequence, we observe 
\[v_n(t)=\sum_{k=0}^n\left(\frac{1}{4}{}^{RL}I^{\alpha}_0-\frac{1}{8}{}^{RL}I^{\alpha+1}_0 \right)^k({}^{RL}I^{\alpha}_0 t^{-\alpha}).\]
By Newton's binomial theorem, the semigroup property of the Riemann-Liouville integral, and passing to the limit as $n\to +\infty$, we have
\begin{align*}
    v(t)&=\sum_{k=0}^{+\infty} \sum_{j=0}^{k}\binom{k}{j}\left( \frac{1}{4}\right)^{k-j}\left( -\frac{1}{8}\right)^{j}{}^{RL}I^{\alpha(k-j)+(\alpha+1)j+\alpha}_0 t^{-\alpha}\\
    &=\Gamma(1-\alpha)\sum_{k=0}^{+\infty} \sum_{j=0}^{k}\binom{k}{j}\left( \frac{1}{4}\right)^{k-j}\left( -\frac{1}{8}\right)^{j}\frac{t^{k\alpha+j}}{\Gamma(k\alpha+j+1)}.
\end{align*}
By the Cauchy product formula with $m=k-j$, it follows that
\begin{align*}
    v(t)&=\Gamma(1-\alpha)\sum_{m=0}^{+\infty}\sum_{j=0}^{+\infty}\frac{(m+j)!}{j!\, m!}\left( \frac{1}{4}\right)^{m}\left( -\frac{1}{8}\right)^{j}\frac{t^{\alpha(m+j)+j}}{\Gamma(\alpha(m+j)+j+1)}\\
    &=\Gamma(1-\alpha)\sum_{m=0}^{+\infty}\sum_{j=0}^{+\infty}\frac{(m+j)!}{j!\, m!}\left( \frac{1}{4}t^\alpha\right)^{m}\left( -\frac{1}{8}t^{\alpha+1}\right)^{j}\frac{1}{\Gamma(\alpha m+(\alpha+1)j+1)}.
\end{align*}
According to \cite [Definition 2.1]{Restrepo}, we arrive at
\[v(t)=\Gamma(1-\alpha)E_{(\alpha,\alpha+1), 1}\left(\frac{1}{4}t^{\alpha},-\frac{1}{8}t^{\alpha+1}\right).\]
On the other hand, applying the Laplace transform to equation \eqref{exam1}, we get
\[
s^{\alpha}\mathcal{L}\{v\}(s)-{}^
{RL}I_0^{1-\alpha}v(0)+\frac{1}{8}s^{-1}\mathcal{L}\{v\}(s)=\frac{1}{4}\mathcal{L}\{v\}(s)+\Gamma(1- \alpha)s^{\alpha-1},\quad s>0.
\]
Using that ${}^{RL}I^{1-\alpha}_0v(0)=0,$ it yields 
\[
s^{\alpha}\mathcal{L}\{v\}(s)+\frac{1}{8}s^{-1}\mathcal{L}\{v\}(s)=\frac{1}{4}\mathcal{L}\{v\}(s)+\Gamma(1- \alpha)s^{\alpha-1}
\]
and 
\begin{align*}
\mathcal{L}\{v\}(s)=\frac{\Gamma(1- \alpha)s^{\alpha-1}}{s^{\alpha}+\frac{1}{8s}-\frac{1}{4}}=\Gamma(1-\alpha)\frac{1}{s(1-\frac{1}{4s^{\alpha}}+\frac{1}{8s^{\alpha+1}})}.
\end{align*}
By \cite [Theorem 2.3]{Restrepo}, we find that
\[v(t)=\Gamma(1-\alpha)E_{(\alpha,\alpha+1), 1}\left(\frac{1}{4}t^{\alpha},-\frac{1}{8}t^{\alpha+1}\right).\]
\end{ex}

\section{Linear integro-differential equations with continuous variable coefficients}\label{c-linear}

In this section, we study the following linear integro-differential equation with continuous variable coefficients:
\begin{equation}\label{C1}
 (\kappa_{0} * v)'(t)
+ \sum_{i=1}^{m} \sigma_{i}(t)\, (\kappa_{i} * v)'(t)
= g(t), \quad t\in[0,T],
\end{equation}
with the initial condition
\begin{equation}\label{C2}
(\kappa_0 * v)(0)=0,
\end{equation}
where $\kappa_0$ is a Sonin kernel, $\sigma_i,g\in C[0,T]$ and $\kappa_i\in AC[0,T]$ for $i=1,\dots,m.$ We look for a continuous solution $v$ to \eqref{C1}-\eqref{C2} in $C^{\kappa_0}$. Note that \eqref{C2} always holds for every $v \in C[0,T].$ 

We start by proving the relationship between conditions $v(0)=0$ and $(\kappa_0*v)(0)=0$ for solutions of \eqref{C1}. Although these conditions are not equivalent for arbitrary functions in $C^{\kappa_0}[0,T]$, they are equivalent for solutions of \eqref{C1}-\eqref{C2}, as shown below.

\begin{prop}\label{propsol}
Let $v\in C^{\kappa_0}$. Then $v$ is a solution of \eqref{C1} with the initial condition $(\kappa_0*v)(0)=0$ if and only if $v$ satisfies \eqref{C1} under the initial condition $v(0)=0.$
\end{prop}
\begin{proof}
Let $v\in C^{\kappa_0}$ be a solution to \eqref{C1} with $v(0)=0.$
From $\kappa_0 \in L^1[0,T]$ and $v\in C[0,T],$ we have 
\begin{align*}
|(\kappa_0*v)(t)|\leqslant\int_0^t|v(s)||\kappa_0(t-s)|\,{\rm d}s\leqslant \|v\|_{\infty}\int_0^t|\kappa_0(r)|\,{\rm d}r.
\end{align*}
As $\displaystyle\lim_{t\to0 ^+}\int_0^t|\kappa_0(u)|\,{\rm d}u =0,$ then $\displaystyle\lim_{t\to0^+}(\kappa_0\ast v)(t)=0.$ Moreover, by Proposition \ref{C3}, we know that $\kappa_0*v\in C[0,T],$ and therefore $(\kappa_0*v)(0)=0.$

Suppose now that $v\in C^{\kappa_0}$ satisfies \eqref{C1} with $(\kappa_0*v)(0)=0.$ From \eqref{C1}, we obtain
\[
(\kappa_0*v)'(t)=g(t)-\sum_{i=1}^{m}\sigma_i(t)(\kappa_i*v)'(t).
\]
Define
\[
h(t):=g(t)-\sum_{i=1}^{m}\sigma_i(t)(\kappa_i*v)'(t).\]
As $v\in C^{\kappa_0}$, one has $\kappa_0*v\in C^1[0,T]$ and $h\in C[0,T].$ By $(\kappa_0*v)(0)=0$, we get
\[
Ih(t)=\int_0^t(\kappa_0*v)'(s)\,{\rm d}s=
(\kappa_0*v)(t)-(\kappa_0*v)(0)=(\kappa_0*v)(t).
\]
From $ (\tilde{\kappa}_0*\kappa_0)(t)=1,$ $(\tilde{\kappa}_0*Ih)(t)=Iv(t).$ Using \cite[Corollary 7.4]{gri}, we arrive at
\[
v(t)=(\tilde{\kappa}_0*Ih)'(t)=Ih(0)\tilde{\kappa}_0(t)+(\tilde{\kappa}_0*h)(t)=(\tilde{\kappa}_0*h)(t).
\]
Since $h\in C[0,T]$ and $\tilde{\kappa}_0 \in L^1[0,T]$, we conclude that $v(0)=(\tilde{\kappa}_0*h)(0)=0.$
\end{proof}
Let us present the principal statement regarding \eqref{C1}-\eqref{C2}. Some lines of the proof are similar to the proof of Theorem \ref{thm-integrable}, so we omit some of the steps.

\begin{thm}\label{The-Continuous}
    If $\sigma_i,g\in C[0,T], \ \kappa_i\in AC[0,T]$ for $i=1,\dots,m,$ then the problem \eqref{C1}-\eqref{C2} has a unique solution $v\in C^{\kappa_0}$, and is explicitly represented as:
    \[v(t)=\sum_{k=0}^{+\infty}(-1)^k \varTheta^k(\tilde{\kappa}_0*g)(t),\]
    where
    \[
    \varTheta w:=\tilde{\kappa}_0 * \left( \sum_{i=1}^m \sigma_i(t)\, (\kappa_i * w)' \right),
    \]
    with $\varTheta^0 := \mathrm{Id}$ and $\varTheta^k$ $(k \geqslant 2)$ is the $k$-fold composition of $\varTheta$.
\end{thm}
\begin{proof}
First, if $v \in C^{\kappa_0}$ satisfies \eqref{C1}-\eqref{C2}, then $u(t)= (\kappa_0 * v)'(t)$ is in $C[0,T]$ and
\begin{equation}\label{C4}
u(t) + \sum_{i=1}^{m} \sigma_{i}(t)(\kappa_i*\tilde{\kappa}_0*u)'(t) = g(t).
\end{equation}
Conversely, from a solution $u\in C[0,T]$ of \eqref{C4}, we obtain a solution $v(t)=(u*\tilde{\kappa}_0)(t)$ in $C^{\kappa_0}$ of \eqref{C1}-\eqref{C2}.

Let $(\lambda_j)_{j \in \mathbb{N}}$ be a sequence of positive integers and define \(T_{\lambda_j}: C[0,T] \to C[0,T]\) by
\begin{equation*}
  T_{\lambda_j}u(t):= e^{-\lambda_j t}\left( g(t)-\sum_{i=1}^m\sigma_i(t)\big(\kappa_i*\tilde{\kappa}_0* u\, e^{\lambda_j (\cdot)}\big)'(t)\right).
\end{equation*}
Obviously, \eqref{C4} is equivalent to  $T_{\lambda_j}\left[e^{-\lambda_j t} u(t)\right] = e^{-\lambda_j t} u(t)$. Let us show that $T_{\lambda_j}$ is a contraction map on $C[0,T]$ for a fixed choice of $\lambda_j,$ for some $j\in\mathbb{N}.$ Let $\lambda_j$ be large enough and fixed. For any $u_1, u_2 \in C[0,T]$ and $t \in [0,T]$, we have 
\begin{align*}
T_{\lambda_j}u_1(t)-T_{\lambda_j}u_2(t)&=-e^{-\lambda_j t}\left(\sum_{i=1}^m\sigma_i(t)\big(\kappa_i*\tilde{\kappa}_0*(u_1-u_2)e^{\lambda_j (\cdot)}\big)'(t)\right).
\end{align*}
Let $w:=u_1-u_2,$ $f:=\tilde{\kappa}_0*we^{\lambda_j (\cdot)}$ and $\kappa_{i,0}=\kappa'_{i}*\tilde\kappa_{0}.$ By \cite[Corollary 7.4]{gri}, we get
\begin{align*}
    (\kappa_i*f)'(t)=\kappa_i(0)(\tilde{\kappa_0}*we^{\lambda_j (\cdot)})(t)+(\kappa_{i,0}*w e^{\lambda_j (\cdot)})(t).
\end{align*}
Hence
\begin{align*}
|T_{\lambda_j}u_1(t) &- T_{\lambda_j}u_2(t)| = e^{-\lambda_j t}\left|\, \sum_{i=1}^m\sigma_i(t)\Big[\kappa_i(0)(\tilde{\kappa}_0*we^{\lambda_j (\cdot)})(t)+(\kappa_{i,0}*we^{\lambda_j (\cdot)})(t)\Big]\,\right|\\  
  &\leqslant \sum_{i=1}^m |\sigma_i(t)|\left[|\kappa_i(0)|\int_0^t|\tilde\kappa_{0}(t-s)|\, |w(s)| e^{-\lambda_j(t-s)}{\rm d}s \right.\\
&\hspace{4cm}\left.+\int_0^t|\kappa_{i,0}(t-s)|\,|w(s)|e^{-\lambda_j(t- s)}{\rm d}s\right]\\
  &\leqslant \|w\|_{\infty}\sum_{i=1}^m\|\sigma_i\|_{\infty}\left[|\kappa_i(0)|\int_0^T|\tilde\kappa_{0}(r)|\, e^{-\lambda_j r }{\rm d}r+\int_0^T|\kappa_{i,0}(r)|\, e^{-\lambda_j r }{\rm d}r\right].
\end{align*}
Young's convolution inequality gives $\kappa_{i,0} \in L^1[0,T]$. By the Lebesgue dominated convergence theorem, it follows that
\[
\int_0^T |\tilde{\kappa}_0(r)| e^{-\lambda_j r}\,{\rm d}r \to 0,
\quad
\int_0^T |\kappa_{i,0}(r)| e^{-\lambda_j r}\,{\rm d}r \to 0,
\quad \text{as } j\to +\infty.
\]
From the previous estimates, we conclude that
\[
|T_{\lambda_j}u_1(t) - T_{\lambda_j}u_2(t)|
\leqslant C(\lambda_j)\,\|u_1-u_2\|_\infty,
\]
where
\[
C(\lambda_j)=\sum_{i=1}^m\|\sigma_i\|_{\infty}
\left(
|\kappa_i(0)|\int_0^T |\tilde{\kappa}_0(r)| e^{-\lambda_j r}\,{\rm d}r
+
\int_0^T |\kappa_{i,0}(r)| e^{-\lambda_j r}\,{\rm d}r
\right).
\]
Since $C(\lambda_j)\to0$ as $j\to +\infty,$ we can find some $0<C(\lambda_j)<1,$ for $j$ sufficiently large, such that
\[
\|T_{\lambda_j}u_1-T_{\lambda_j}u_2\|_\infty
\leqslant C(\lambda_j)\|u_1-u_2\|_\infty.
\]
By the Banach fixed-point theorem, there exists a unique solution $u\in C[0,T]$ to \eqref{C4} which can be approximated by the following sequence
    \begin{align*}
        u_0(t)&=g(t), \nonumber  \\
        u_n(t)&= g(t) - \sum _{i=1}^m\sigma_i(t)(\kappa_{i}*\tilde{\kappa}_0*u_{n-1})'(t), \quad n \geqslant 1. %\label{C6}
    \end{align*}
Set $v_n(t) = (\tilde{\kappa}_0 * u_n)(t)$ for all $n$, then 
    \begin{align}
        v_0(t) &= (\tilde{\kappa}_0*g)(t), \nonumber \\
        v_n(t)&= v_0(t)- \tilde \kappa_0*\sum _{i=1}^m\sigma_i(t)(\kappa_i*v_{n-1})'(t), \quad n\geqslant 1.  {\label{C8}}
    \end{align} 
It can be proved that $\{v_n\}_{n\geqslant 0}$ converges in $C^{\kappa_0}$, and the limit function $v=\tilde{\kappa}_0*u$ is the unique solution of \eqref{C1} in $C^{\kappa_0}.$ We can construct the explicit solution of $v$ using recurrence \eqref{C8} as was done for \eqref{Suc1}. The solution $v$ is given by $$v(t)=\lim_{n\to\infty}v_n(t)=\sum_{k=0}^{+\infty}(-1)^k \varTheta^k(\tilde{\kappa}_0*g)(t),$$
where the limit is taken in $\|\cdot\|_{C^{\kappa_0}}$ and the series is convergent in this norm.
\end{proof}
\begin{rem}[Constant coefficients]
When the coefficients $\sigma_i(t)$ in \eqref{C1} are constant,
we do not obtain weaker conditions on the kernels or a simpler
formula for the solution. Since this case is already considered in Theorem 5.2, we will not analyze it separately.
\end{rem}
\begin{rem}
Let us discuss the following question: Can problem \eqref{C1} admit a solution $v\in C^{\kappa_0}$ with an initial condition $v(0)=a\neq 0\in\mathbb{R}$?

Using our approach, this is not possible. On the contrary, suppose that there exists $v\in C^{\kappa_0}$ satisfying \eqref{C1} with $v(0)=a\neq0.$ Define $u(t)=({\kappa}_0*v)'(t).$ We know that $u \in C[0,T]$ and $Iu(t) = (\kappa_0 * v)(t)- (\kappa_0 * v)(0).$ Since $v\in C[0,T]$ and ${\kappa}_0 \in L^1[0,T]$, it follows that $({\kappa}_0*v)(0)=0.$ From $(\tilde{\kappa}_0*\kappa_0)(t)=1,$ we have $(\tilde{\kappa}_0*Iu)(t)=Iv(t).$ By \cite[Corollary 7.4]{gri}, we obtain
\[v(t)=Iu(0)\tilde{\kappa}_0(t)+(\tilde{\kappa}_0*u)(t)=(\tilde{\kappa}_0*u)(t).\]
Therefore, $v(t)=(\tilde{\kappa}_0*u)(t)$ and $v(0)=0.$ This is a contradiction.

The latter brief exposition shows that our model is restricted for initial conditions to be zero. This limitation follows from the spaces considered $C^{\kappa_0}.$ This opens a new site of investigation on the proper functional spaces to study integro-differential equations with continuous solutions.      

\end{rem}

The following example is a direct application of Theorem \ref{The-Continuous}.

\begin{ex}
Take $\kappa_0(t)=t^{-\alpha}/\Gamma(1-\alpha)$ with $\alpha\in(0,1),$ $\sigma_1(t)=1$, $\kappa_1(t)=1$ and $g(t)=e^{\lambda t}$ for some $\lambda>0.$ With the latter functions, equation \eqref{C1} turns to
\begin{equation}\label{ex_2}
   \left(\frac{t^{-\alpha}}{\Gamma(1-\alpha)}*v\right)'(t)+(1*v)'(t)=e^{\lambda t}, \quad t\in [0,2],      
    \end{equation}
 under the initial condition 
 \begin{equation}\label{ex_2-c}
  \left(\frac{t^{-\alpha}}{\Gamma(1-\alpha)}*v\right)(0)=0.   
 \end{equation}
By Theorem \ref{The-Continuous}, \eqref{ex_2}-\eqref{ex_2-c} have a unique solution $v$ in $C^{\frac{t^{-\alpha}}{\Gamma(1-\alpha)}}$, that is given by
    \[ v(t)=\sum_{k=0}^{+\infty}(-1)^k \varTheta^k\left(\frac{t^{\alpha-1}}{\Gamma(\alpha)}*e^{\lambda t}\right)=\sum_{k=0}^{+\infty}(-1)^k \varTheta^k\left({}^{RL}I^{\alpha}_0e^{\lambda t}\right).
    \]
For $w={}^{RL}I^{\alpha}_0e^{\lambda t},$ one has
\[ 
\varTheta w = \tilde{\kappa}_0 * \left( \sum_{i=1}^{1} \sigma_i(t) (\kappa_i * w)' \right) = \frac{t^{\alpha-1}}{\Gamma(\alpha)} * \left( 1 * w \right)'=\left(\frac{t^{\alpha-1}}{\Gamma(\alpha)}*w\right)={}^{RL}I^{\alpha}_0w.\]
By the semigroup property of the Riemann-Liouville fractional integral, the $k$-fold composition of the operator takes the form $\varTheta^k w= {}^{RL}I^{k \alpha}_0w.$ Substituting this representation of $\Theta^k$ into the series solution, we obtain
\begin{align*}
v(t) &= \sum_{k=0}^{+\infty} (-1)^k  \, {}^{RL}I_{0}^{k\alpha} \big({}^{RL}I_{0}^{\alpha}e^{\lambda t}\big) = \sum_{k=0}^{+\infty} (-1)^k\, {}^{RL}I_{0}^{\alpha(k+1)} e^{\lambda t}.
\end{align*}
We now compute 
\[ {}^{RL}I_{0}^{\alpha(k+1)} e^{\lambda t}= \frac{\lambda^{-\alpha(k+1)}e^{\lambda t}}{\Gamma(\alpha(k+1))}\int_0^{\lambda t} r^{\alpha(k+1)-1} e^{-r}{\rm d}r=\frac{\lambda^{-\alpha(k+1)}e^{\lambda t}}{\Gamma(\alpha(k+1))}\gamma(\alpha(k+1),\lambda t),
\]
where $\gamma(\cdot,\cdot)$ is the lower incomplete gamma function \cite[Formula~1.5.27]{kilbas}. Hence
\[ v(t)=\sum_{k=0}^{+\infty} (-1)^k \frac{\lambda^{-\alpha(k+1)}e^{\lambda t}}{\Gamma(\alpha(k+1))}\gamma(\alpha(k+1),\lambda t).\]
Observe that directly from the representation of $v,$ we have $v(0)=0$ and 
\begin{align*}
    {}^{RL}I_{0}^{1-\alpha} v(t)&=\sum_{k=0}^{+\infty}(-1)^k\, {}^{RL}I_{0}^{1-\alpha} ( {}^{RL}I_{0}^{\alpha(k+1)} e^{\lambda t})=\sum_{k=0}^{+\infty}(-1)^k{}^{RL}I_{0}^{\alpha k+1} e^{\lambda t} \\&=\sum_{k=0}^{+\infty}(-1)^k\frac{\lambda^{-(\alpha k+1)}e^{\lambda t}}{\Gamma(\alpha k+1)}\gamma(\alpha k+1,\lambda t),
    \end{align*}
which implies ${}^{RL}I_{0}^{1-\alpha} v(0)=0.$ 

Again, for completeness, let us find the solution by applying the Laplace transform to equation \eqref{ex_2}. Thus
\[
s^{\alpha}\mathcal{L}\{v\}(s)-{}^{RL}I_0^{1-\alpha}v(0)+\mathcal{L}\{v\}(s)=\frac{1}{s-\lambda},\quad s>0,
\]
and
\[\mathcal{L}\{v\}(s)=\frac{1}{s^\alpha+1}\frac{1}{s-\lambda}.\]
By properties of the Laplace transform, we get
\[v(t)=\mathcal{L}^{-1}\left\{\frac{1}{s^\alpha+1}\right\}(t)*\mathcal{L}^{-1}\left\{\frac{1}{s-\lambda}\right\}(t).\]
By \cite [Theorem 2.3]{Restrepo}, 
\[
v(t)=\left(t^{\alpha-1}E_{\alpha,\alpha}(-t^{\alpha})*e^{\lambda t}\right)(t)=\sum_{k=0}^{+\infty} (-1)^k \, {}^{RL}I_0^{\alpha(k+1)}e^{\lambda t}.
\]
Here, the representation obtained through the Laplace transform coincides with the solution derived from Theorem \ref{The-Continuous}. Problem \eqref{ex_2}-\eqref{ex_2-c} can be also rewritten with the Riemann-Liouville integro-differential operators as
\[
{}^{RL}D_0^{\alpha}v(t)+v(t)=e^{\lambda t},\quad {}^{RL}I_{0}^{1-\alpha} v(0)=0.
\]
\end{ex}

\begin{ex}[A fractional kinetic model]
\label{ex:fractional-kinetics}
Haubold and Mathai considered the following fractional kinetic model \cite[Equation~21]{HauboldMathai2000}:
\begin{equation}\label{eq:fractional-kinetics}
N(t)-N_0=-c^\alpha\,{}^{RL}I_0^\alpha N(t),
\qquad 0<t<T,\quad 0<\alpha<1,
\end{equation}
where $c>0$ is the destruction rate coefficient, $N(t)$ is the number density of the species, and $N_0=N(0)$ is its initial value.

Assume that $N\in C[0,T]$. Applying the Riemann--Liouville fractional derivative of order $\alpha$ to both sides of \eqref{eq:fractional-kinetics}, we obtain
\[
{}^{RL}D_0^\alpha(N-N_0)(t)=-c^\alpha N(t).
\]
Setting $v(t)=N_0-N(t)$, we have $v(0)=0$; and the previous equation becomes
\[
{}^{RL}D_0^\alpha v(t)+c^\alpha v(t)
=c^\alpha N_0.
\]
We can rewrite this equation as 
\[
\left(
\frac{t^{-\alpha}}{\Gamma(1-\alpha)}*v
\right)'(t)
+c^\alpha(1*v)'(t)=c^\alpha N_0.
\]
Here, $\kappa_0(t)=t^{-\alpha}/\Gamma(1-\alpha)$, $\tilde{\kappa}_0(t)=t^{\alpha-1}/\Gamma(\alpha),$ $\kappa_1(t)=1$, $\sigma_1(t)=c^\alpha$ and $g(t)=c^\alpha N_0$. Theorem~\ref{The-Continuous} guarantees the existence of a unique solution $v\in C^{\kappa_0}$ of the above equation. In this case, the operator $\varTheta$ takes the form
\[
\varTheta w=\frac{t^{\alpha-1}}{\Gamma(\alpha)}*\left(c^\alpha(1*w)'\right)=c^\alpha\,{}^{RL}I_0^\alpha w.
\]
Using the semigroup property of the Riemann-Liouville fractional integral:
\[
\varTheta^k w
=c^{k\alpha}\,{}^{RL}I_0^{k\alpha}w.
\]
Substituting these expressions into the series
representation of the solution from Theorem~\ref{The-Continuous}, we have
\begin{align*}
v(t)
&=\sum_{k=0}^{+\infty}(-1)^k c^{k\alpha}{}^{RL}I_0^{k\alpha}\left(\frac{t^{\alpha-1}}{\Gamma(\alpha)}*c^\alpha N_0\right)(t)=\sum_{k=0}^{+\infty}(-1)^k c^{(k+1)\alpha}\left({}^{RL}I_0^{(k+1)\alpha}N_0\right)(t)\\
&=N_0\sum_{k=0}^{+\infty}(-1)^k\frac{(ct)^{(k+1)\alpha}}{\Gamma((k+1)\alpha+1)}
=-N_0\sum_{j=1}^{+\infty}\frac{(-c^\alpha t^\alpha)^j}{\Gamma(j\alpha+1)}
=N_0\left(1-E_\alpha(-c^\alpha t^\alpha)\right).
\end{align*}
Consequently, recalling that $N(t)=N_0-v(t)$, we have $N(t)=N_0E_\alpha(-c^\alpha t^\alpha).$
\end{ex}

\section{Nonlinear integro-differential equations with continuous variable coefficients} \label{nonlinear-continuous}

We generalize equation \eqref{C1} by incorporating a nonlinear term. We consider the following equation:
\begin{equation}\label{CN1}
 (\kappa_0 * v)'(t)
+ \sum_{i=1}^{m} \sigma_i(t)\,(\kappa_i * v)'(t)
= g(t,v(t)), \quad t \in [0,T],
\end{equation} 
with the initial condition 
\begin{equation}\label{CN2}
    (\kappa_0 * v)(0)=0.
\end{equation}
We suppose that $\kappa_0$ is a Sonin kernel, $\sigma_i \in C[0,T]$ and $\kappa_i \in AC[0,T]$ for $i=1,\dots,m.$ Let $\Delta\subseteq\mathbb{R}$ be an open set and  $g:[0,T]\times\Delta\to\mathbb{R}$ be a continuous function such that 
\begin{equation}\label{CN3}
|g(t,v_1(t))-g(t,v_2(t))| \leqslant \lambda|v_1(t)-v_2(t)|, \quad  t\in[0,T], \quad v_1(t),v_2(t)\in\Delta,
\end{equation}
for some constant $\lambda>0.$ 
\begin{rem}
Observe that similar to Proposition \ref{propsol}, for $v\in C^{\kappa_0},$ it can be easily seen that $v$ is a solution of \eqref{CN1}-\eqref{CN2} is a necessary and sufficient condition for $v$ to be a solution of \eqref{CN1} with $v(0)=0.$ Therefore, equation \eqref{CN1} can be solved using either of these two initial conditions.
\end{rem}

Our goal is to prove that under \eqref{CN2}–\eqref{CN3}, \eqref{CN1} has a unique solution $v \in C^{\kappa_0}.$ The following proof is a continuation of the one given in Theorem \ref{The-Continuous}. 

 \begin{thm}\label{non-linear-thm}
       Assume that $\sigma_i \in C[0,T]$, $\kappa_i \in AC[0,T]$, for $i=1,\dots,m$ and $g$ is continuous on $[0,T]\times\Delta$. If conditions \eqref{CN2}-\eqref{CN3} are satisfied, then equation \eqref{CN1} admits a unique solution $v \in C^{\kappa_0}$.
 \end{thm}
 
\begin{proof}
Observe that if $v \in C^{\kappa_0}$ is a solution of \eqref{CN1}, then $u(t) = (\kappa_0 * v)'(t)\in C[0,T]$ satisfies the following equation 
\begin{equation}\label{CN4}
u(t) + \sum_{i=1}^{m} \sigma_{i}(t) (\kappa_i*\tilde{\kappa}_0*u)'(t) = g\left(t,(\tilde\kappa_0*u)(t)\right),
\end{equation}
and vice versa. This means, if $u \in C[0,T]$ is a solution to \eqref{CN4}, then $v(t)=(\tilde{\kappa}_0*u)(t)$ is in $C^{\kappa_0}$ and satisfies \eqref{CN1}.

Consider a sequence $(\alpha_j)_{j \in \mathbb{N}}$ of positive integers and $\widetilde{T}_{\alpha_j}: C[0,T] \to C[0,T]$ as
\begin{equation*}
  \widetilde{T}_{\alpha_j}u(t):= e^{-\alpha_j t}\left( g(t,(\tilde{\kappa}_0*(ue^{\alpha_j(\cdot)}))(t))-\sum_{i=1}^m\sigma_i(t)(\kappa_i*\tilde{\kappa}_0* (ue^{\alpha_j(\cdot)}))'(t)\right),
\end{equation*}
where $(ue^{\alpha_j(\cdot)})(t)=u(t)e^{\alpha_jt}.$ By \eqref{CN4}, we know that $\widetilde{T}_{\alpha_j}\big[e^{-\alpha_j t}u(t)\big]=e^{-\alpha_j t}u(t)$. Let us check that for a suitable choice of $\alpha_j$, $\widetilde{T}_{\alpha_j}$ is a contraction on $C[0,T].$ For $u_1, u_2 \in C[0,T]$ we have
\begin{align*}
\widetilde{T}_{\alpha_j}u_1(t)-\widetilde{T}_{\alpha_j}u_2(t)&=e^{-\alpha_jt}\left[g(t,(\tilde{\kappa}_0*(u_1e^{\alpha_j(\cdot)}))(t))-g(t,(\tilde{\kappa}_0*(u_2e^{\alpha_j(\cdot)}))(t))\right] \\
&\hspace{1cm}-e^{-\alpha_jt}\left[\sum_{i=1}^m\sigma_i(t)\big[(\kappa_i*\tilde{\kappa}_0*((u_1-u_2)e^{\alpha_j(\cdot)}))(t)\big]'\right].
\end{align*}
Let $w=u_1-u_2$ and $f=\tilde{\kappa}_0*(we^{\alpha_j(\cdot)})$. By \cite[Corollary 7.4]{gri} and \eqref{CN3}, we find that
\begin{align*}
|\widetilde{T}_{\alpha_j}&u_1(t) - \widetilde{T}_{\alpha_j}u_2(t)| \leqslant    e^{-\alpha_j t}\left|g(t,(\tilde{\kappa}_0*(u_1e^{\alpha_j(\cdot)}))(t))-g(t,(\tilde{\kappa}_0*(u_2e^{\alpha_j(\cdot)}))(t))\right|\\
 &\hspace{1cm}+e^{-\alpha_j t} \left|\sum_{i=1}^m\sigma_i(t)\left[\kappa_i(0)(\tilde{\kappa}_0*(we^{\alpha_j (\cdot)}))(t)+(\kappa_{i,0}*(we^{\alpha_j (\cdot)}))(t)\right]\right| \\
 &\leqslant e^{-\alpha_j t}\lambda|(\tilde{\kappa}_0*(u_1e^{\alpha_j(\cdot)}))(t)-(\tilde{\kappa}_0*(u_2e^{\alpha_j(\cdot)}))(t)|\\
&\hspace{1cm}+\sum_{i=1}^m|\sigma_i(t)|\left[e^{-\alpha_j t}\left|\kappa_i(0)(\tilde{\kappa}_0*(we^{\alpha_j (\cdot)}))(t)\right|+e^{-\alpha_j t}\left|(\kappa_{i,0}*(we^{\alpha_j (\cdot)}))(t)\right|\right]\\
  &\leqslant \lambda\int_0^t |\tilde{\kappa}_0(t-s)||w(s)|\,e^{-\alpha_j(t-s)}{\rm d}s\\
  &\hspace{1cm}+\sum_{i=1}^m |\sigma_i(t)|\Big[|\kappa_i(0)|\int_0^t|\tilde\kappa_{0}(t-s)|\, |w(s)| e^{-\alpha_j(t-s)}{\rm d}s\\
  &\hspace{2cm}+\int_0^t|\kappa_{i,0}(t-s)|\,|w(s)|e^{-\alpha_j(t- s)}\,{\rm d}s\Big]\\
   &\leqslant \|w\|_{\infty}\Bigg[\lambda\int _0^T|\tilde{\kappa}_0(r)|e^{-\alpha_j r}{\rm d}r+\sum_{i=1}^m\|\sigma_i\|_{\infty}\Big[|\kappa_i(0)|\int_0^T|\tilde\kappa_{0}(r)|\, e^{-\alpha_j r }{\rm d}r \\
  &\hspace{1cm}+\int_0^T|\kappa_{i,0}(r)|\, e^{-\alpha_jr }{\rm d}r\Bigg].
 \end{align*}
The Lebesgue dominated convergence theorem gives
\begin{equation*}
\int_0^T |\tilde{\kappa}_0(r)|e^{-\alpha_j r}\,{\rm d}r \to 0,
\qquad
\int_0^T |\kappa_{i,0}(r)|e^{-\alpha_j r}\,{\rm d}r \to 0, \quad \text{as}\quad j\to+\infty.
\end{equation*}
Hence $|\widetilde{T}_{\alpha_j}u_1(t) - \widetilde{T}_{\alpha_j}u_2(t)| \leqslant  C(\alpha_j) \,\|u_1-u_2\|_{\infty},$ where 
\begin{align*}
C(\alpha_j)&=\lambda\int_0^T|\tilde{\kappa}_0(r)|e^{-\alpha_jr} {\rm d}r \\
&\hspace{1cm}+\sum_{i=1}^m\|\sigma_i\|_{\infty}\left[ |\kappa_i(0)|\int_0^T|\tilde{\kappa}_0(r)|e^{-\alpha_jr}{\rm d}r+\int_0^T|\kappa_{i,0}(r)|e^{-\alpha_jr}{\rm d}r\right].
\end{align*}
Since the above integrals converge to zero as $j \to +\infty$, we have $C(\alpha_j)\to 0$ as $j \to +\infty$. For $j$ sufficiently large, we obtain
\[ \|\widetilde{T}_{\alpha_j}u_1 -\widetilde{T}_{\alpha_j}u_2\|_{\infty} \leqslant C(\alpha_j)\|u_1 - u_2\|_{\infty},\quad\text{for some}\quad 0 < C(\alpha_j) < 1.\]
By the Banach fixed point theorem, $\tilde T_{\alpha_j}$ admits a unique fixed point $u\in C[0,T].$ So, problem \eqref{CN1}-\eqref{CN2}-\eqref{CN3} have a unique solution $v\in C^{\kappa_0}$.
\end{proof}

\begin{rem}\label{nolinear-rem}
From Theorem \ref{non-linear-thm}, the solution $u$ of \eqref{CN4} can be obtained as the limit of the following iterative sequence
\begin{align*}%\label{CN5}
        u_0(t)&=g(t,0) ,\nonumber \\
        u_n(t)&= g\big(t,(\tilde{\kappa}_0*u_{n-1})(t)\big) - \sum _{i=1}^m\sigma_i(t)(\kappa_{i}*\tilde{\kappa}_0*u_{n-1})'(t), \quad n\geqslant1,
\end{align*}
which converges in $\|\cdot\|_{\infty}.$ Let $v_n=(\tilde{\kappa}_0*u_n)(t)$ for every $n$, then
    \begin{align}\label{CN7}
        v_0(t) &= (\tilde{\kappa}_0*g(\cdot,0))(t), \nonumber \\
        v_n(t)&= (\tilde{\kappa}_0*g(\cdot,v_{n-1}(\cdot)))(t)- \tilde\kappa_0*\sum _{i=1}^m\sigma_i(t)(\kappa_i*v_{n-1})'(t).
    \end{align} 
It is easy to see that $v_n \in C^{\kappa_0}$ for all $n.$ Theorem \ref{The-Continuous} shows that the sequence $\{v_n\}_{n\geqslant0}$ converges in $C^{\kappa_0}$. Consequently, its limit $v\in C^{\kappa_0}$ is the unique solution of \eqref{CN1}.
\end{rem}

In the next example, we use Remark \ref{nolinear-rem} to find the solution iteratively.  

 \begin{ex}
Assume $\kappa_0(t)=t^{-\alpha}/\Gamma(1-\alpha)$ for $\alpha\in(0,1),$  $\sigma_1(t)=1$, $\kappa_1(t)=t$ and $g(t,v(t))= v(t)+\cos(t).$ Hence, equation \eqref{CN1} turns to
\begin{equation}\label{ex3}
\left( \ \frac{t^{-\alpha}}{\Gamma(1-\alpha)}  * v \right)'(t)+\left(t*v\right)'(t)= v(t)+\cos(t), \quad t \in [0,10],
\end{equation}
with initial condition $\left(\frac{t^{-\alpha}}{\Gamma(1-\alpha)}*v\right)(0)=0.$ Clearly $|v_1(t)+\cos(t)-v_2(t)-\cos(t)|=|v_1(t)-v_2(t)|.$ By Theorem \ref{non-linear-thm}, \eqref{ex3} has a unique solution $ v \in C^{\frac{t^{-\alpha}}{\Gamma( 1-\alpha)}}.$ From \eqref{CN7}, we get the following iterative sequence
\begin{align*}
    v_0(t)&=\left(\frac{t^{\alpha-1}}{\Gamma(\alpha)}*\cos(t)\right)(t),\\
    v_{n}(t)&=\left(\frac{t^{\alpha-1}}{\Gamma(\alpha)}*( v_{n-1}+\cos(t))\right)(t)-\left(\frac{t^{\alpha-1}}{\Gamma(\alpha)}*(t*v_{n-1})'\right)(t),
\end{align*}
or which is the same
\begin{align*}
    v_0(t)&={}^{RL}I^{\alpha}_0 \cos(t),\\
    v_n(t)&={}^{RL}I^{\alpha}_0 v_{n-1}(t)+{}^{RL}I^{\alpha}_0 \cos(t)-{}^{RL}I^{\alpha+1}_0 v_{n-1}(t).
\end{align*}
Successive iterations of the previous sequence give
\[v_n(t)=\sum_{k=0}^n\left({}^{RL}I^{\alpha}_0-{}^{RL}I^{\alpha+1}_0 \right)^k\big({}^{RL}I^{\alpha}_0 \cos(t)\big).\]
By Newton's binomial theorem, \cite[Theorem 2.5]{samko} and letting $n\to+\infty$, we have
\begin{align*}
v(t)&=\sum_{k=0}^{+\infty} \sum_{j=0}^{k}\binom{k}{j}( -1)^{j}\,{}^{RL}I^{\alpha(k-j)+(\alpha+1)j+\alpha}_0 \cos(t).
\end{align*}
By the Cauchy product formula with $m=k-j$, we obtain
\begin{align*}
v(t)&=\sum_{m=0}^{+\infty}\sum_{j=0}^{+\infty}( -1)^{j}\frac{(m+j)!}{j!\, m!}{}^{RL}I^{m\alpha+(\alpha+1)j+\alpha}_0\cos(t)\\
    &=\int_0^t\cos(t-s)\sum_{m=0}^{+\infty}\sum_{j=0}^{+\infty}(-1)^j\frac{(m+j)!}{j!\, m!}\frac{s^{m\alpha+(\alpha+1)j+\alpha-1}}{\Gamma(m\alpha+(\alpha+1)j+\alpha)}\,{\rm d}s\\
    &=\int_0^t\cos(t-s)s^{\alpha-1}\sum_{m=0}^{+\infty}\sum_{j=0}^{+\infty}\frac{(m+j)!}{j!\, m!}\frac{(s^\alpha)^m (-s^{\alpha+1})^{j}}{\Gamma(m\alpha+(\alpha+1)j+\alpha)}\,{\rm d}s.
\end{align*}
From \cite [Definition 2.1]{Restrepo}, we find that
\[v(t)=\int_0^t\cos(t-s)s^{\alpha-1}E_{(\alpha,\alpha+1),\alpha}(s^\alpha,-s^{\alpha+1})\,{\rm d}s.\]
Observe that \eqref{ex3} can be rewritten as
\begin{equation}\label{ex3.1}
   {}^{RL}D^\alpha_0v(t)+ I v(t)= v(t)+\cos(t).
\end{equation}
Also, by applying the Laplace transform to \eqref{ex3.1} yields
\[s^{\alpha}\mathcal{L}\{v\}(s)-{}^
{RL}I_0^{1-\alpha}v(0)+s^{-1}\mathcal{L}\{v\}(s)=\mathcal{L}\{v\}(s)+\frac{s}{s^2+1},\quad s>0.
\]
We know that ${}^{RL}I^{1-\alpha}_0v(0)=0$ and hence
\[
s^{\alpha}\mathcal{L}\{v\}(s)+s^{-1}\mathcal{L}\{v\}(s)=\mathcal{L}\{v\}(s)+\frac{s}{s^2+1},
\]
and 
\[
\mathcal{L}\{v\}(s)=\frac{s^2}{(s^2+1)(s^{\alpha+1}-s+1)}.
\]
Therefore
\[{v}(t)=\mathcal{L}^{-1}\left\{\frac{s^2}{(s^2+1)(s^{\alpha+1}-s+1)}\right\}=\mathcal{L}^{-1}\left\{ \frac{s}{s^2+1}\right\}*\mathcal{L}^{-1}\left\{ \frac{s}{s^{\alpha+1}-s+1}\right\}.\]
By  \cite [Theorem 2.3]{Restrepo}, we get
\[v(t)=\int_0^t\cos(t-s)s^{\alpha-1}E_{(\alpha,\alpha+1), \alpha}(s^{\alpha},-s^{\alpha+1})\,{\rm d}s.\]
\end{ex}

\section{Conclusions}

This work studied a class of linear and nonlinear integro-differential equations with Sonin and absolutely continuous kernels, considering both integrable and continuous variable coefficients. The main objectives were to establish the minimum conditions for the solvability of these problems. In the linear case, we prove the existence and uniqueness of the solution and obtain an analytical representation of it. In the nonlinear case, we guaranty the existence and uniqueness of the solution, as well as its approximation, using an iterative process.

Our viewpoint sheds some new light on the equivalence of conditions $v(0)=0$ and $(\kappa_0*v)(0)=0$ for solutions of \eqref{C1}, see Proposition \ref{propsol}. This property allows us to reformulate the original problem using a more natural initial condition, expressed directly in terms of the unknown function. In the same context, we do not require that $k_i(0)=0$ for $i=1,\ldots,m.$

\medskip The results obtained in this work open up different lines of future research:
\begin{enumerate}

\item A first natural problem is to weaken the assumptions that guaranty the contraction of the operators involved in the existence and uniqueness results for the linear and nonlinear equations.

\item Another problem is the study of equations under nonhomogeneous initial conditions. In particular, it would be of interest to analyze problems with conditions of the form $(\kappa_0*v)(0)\neq 0 $ or $v(0)\neq 0,$ which in some cases are not considered in this paper.

\item Consider some other techniques and tools to study integro-differential equations with Sonin kernels and variable coefficients. Here, we just consider the Banach fixed point theorem and the method of successive approximations.   

\end{enumerate}

\section{Data availability statement and acknowledgments}
Data sharing not applicable to this article as no datasets were generated or analyzed during the current study. No new data was collected or generated during the course of this research.


\begin{thebibliography}{99}

\bibitem{Restrepo}{A. Abilassan, J.E. Restrepo, D. Suragan. On a variant of multivariate Mittag-Leffler's function arising in the Laplace transform method. Integr. Transf. Spec. F., 34(3), (2022), 244—260.}

\bibitem{kober} F. Al-Musalhi, A. Fernandez. Fractional differential equations involving Erd\'elyi-Kober derivatives with variable coefficients. Fract. Calc. Appl. Anal., 28, (2025), 1252--1281 .

\bibitem{transmuted} F. Al-Musalhi. Solving fractional differential equations with variable coefficients using transmutation relations. Comp. Appl. Math., 45(26), (2026), 1--21.

\bibitem{vcserbia1} T.M. Atanackovi\'c, B. Stankovi\'c. Linear fractional differential equation with variable coefficients I. Bull. de l Acad. Serbe Sci. Arts, Cl. Math., 38, (2013), 27--42.

\bibitem{vcserbia2} T.M. Atanackovi\'c, B. Stankovi\'c. Linear fractional differential equation with variable coefficients II. Bull. de l Acad. Serbe Sci. Arts, Cl. Math., 39, (2014), 53--78.

\bibitem{diethelm}
K. Diethelm. The analysis of fractional differential equations: An application-oriented exposition using differential operators of Caputo type. Springer, Heidelberg, 2010.

\bibitem{first} M.M. Dzhrbashyan, A.B. Nersessyan. Fractional derivatives and Cauchy problem for differential equations of fractional order. Izv. AN Arm. SSR. Mat,. 3, (1968).

\bibitem{faustino} N. Faustino. On fractional semidiscrete Dirac operators of L\'e-Leblond type. Math. Nachr., 296(7), (2023), 2758--2779.

\bibitem{FRS:AB} A. Fernandez, J.E. Restrepo, D. Suragan. Linear differential equations with variable coefficients and Mittag-Leffler kernels. Alexandria Eng. J., 61, (2022), 4757--4763.


\bibitem{FRS} A. Fernandez, J.E. Restrepo, D. Suragan. On linear fractional differential equations with variable coefficients. Appl. Math. Comput., 432, (2022), \#127370.

\bibitem{FRS:P} A. Fernandez, J.E. Restrepo, D. Suragan. Prabhakar-type linear differential equations with variable coefficients. Differ. Integral Equ., 35(9/10), (2022), 581--610.

\bibitem{garra} R. Garra. Analytic solution of a class of fractional differential equations with variable coefficients by operatorial methods. Commun. Nonlinear Sci. Numer. Simul., 17(4), (2012),1549--1554.

\bibitem{GarrappaMainardiMaione2016} R. Garrappa, F. Mainardi, G. Maione. Models of dielectric relaxation based on completely monotone functions.
Fract. Calc. Appl. Anal. 19(5), (2016), 1105--1160.

\bibitem{Gorenflo}R. Gorenflo, A. A. Kilbas, F. Mainardi, S.V. Rogosin. Mittag-Leffler Functions, Related Topics and Applications, 2nd ed. Springer Monographs in Mathematics, Springer, New York, 2020.

\bibitem{gri} G. Gripenberg, S.O. Londen, O.J. Staffans. Volterra Integral and Functional Equations. Cambridge University Press, Cambridge, 1990.


\bibitem{sonine-arran} Y. G\"uzoglu Y\"ucel, A. Fernandez, N.I. Mahmudov. Fractional differential equations with continuous variable coefficients and Sonin kernels. Commun. Nonlinear Sci. Numer. Simul., 152, (2026), \#109136.

\bibitem{Koeller1984} R.C. Koeller. Applications of fractional calculus to the theory of viscoelasticity. J. Appl. Mech. 51(2), (1984), 299--307.

\bibitem{Luchko2021} Y.~Luchko. General fractional integrals and derivatives with the Sonine kernels, Mathematics, 9(6):594, (2021), 1--17.

\bibitem{Yu-sym} Y.~Luchko. General fractional integrals and derivatives of arbitrary order. Symmetry, 13(5), (2021), $\#$ 755.

\bibitem{Yu-fcaa} Y.~Luchko. Operational calculus for the general fractional derivatives with the Sonine kernels. Fract. Calc. Appl. Anal., 24(2), (2021), 338--375. 

\bibitem{Yu-1} Y. Luchko. Convolution series and the generalized convolution Taylor formula. Fract. Calc. Appl. Anal., 25, (2022), 207--228.
 
\bibitem{Yu-2} Y. Luchko. General fractional integrals and derivatives and their applications. Phys D, 455, (2023), $\#$133906.

\bibitem{HauboldMathai2000}
H.J. Haubold, A.M. Mathai. The fractional kinetic equation and thermonuclear functions. Astrophys. Space Sci., 273, (2000), 53--63.

\bibitem{hilfer}
R. Hilfer, ed. Applications of Fractional Calculus in Physics. World Scientific, Singapore, 2000.

\bibitem{kilbasalpha} A.A. Kilbas, M. Rivero, L. Rodriguez-Germa, J.J. Trujillo. $\alpha$-Analytic solutions of some linear fractional differential equations with variable coefficients, Appl. Math. Comput. 187 (2007), 239--249.

\bibitem{kilbas} A.A. Kilbas, H.M. Srivastava, J.J. Trujillo. Theory and Applications of Fractional Differential Equations. North-Holland Mathematics Studies, vol. 204, Amsterdam, 2006.

\bibitem{RL} M. Kim, H.-C. O. Explicit representations of Green's function for linear fractional differential operator with variable coefficients. J. Fract. Calc. Appl. 5(1) (2014), 26--36.

\bibitem{Kochubei} A.N. Kochubei. General fractional calculus, evolution equations, and renewal processes. Integr. Equ. Oper. Theory, 71, (2011), 583--600.

\bibitem{hilfer-psi} F. Li, H. Wang. $\psi-$Hilfer type linear fractional differential equations with variable coefficients. Fract. Calc. Appl. Anal., 28, (2025), 807--838.

\bibitem{Li-nonlinear} F. Li, H. Wang. Constructions of analytical solutions to nonlinear fractional differential equations with variable coefficients. Chaos Soliton Fract., 205, (2026), \#117873. 

\bibitem{mainardi} F. Mainardi. Fractional Calculus and Waves in Linear Viscoelasticity. Imperial College Press, London, 2010.

\bibitem{MainardiSpada2011}
F. Mainardi, G. Spada. Creep, relaxation and viscosity properties for basic fractional models in rheology. Eur. Phys. J. Spec. Top. 193(1), (2011), 133--160.

%\bibitem{Nielsen1997} O.A. Nielsen. An Introduction to Integration and Measure Theory. Wiley--Interscience, 1997.

%\bibitem{miller} K.S. Miller, B. Ross. An Introduction to the Fractional Calculus and Fractional Differential Equations. John Wiley, New York, 1993.


\bibitem{oldham} K.B. Oldham, J. Spanier. The Fractional Calculus. Academic Press, New York, 1974.

\bibitem{greenfunction} S. Pak, H. Choi, K. Sin. Representation of solutions of linear inhomogeneous Caputo fractional differential equation with continuous variable coefficient by Green function. Fractional Differ. Calc., 7(2), (2017), 357--369.

\bibitem{analitical} S. Pak, H. Choi, K. Sin, K. Ri. Analytical solutions of linear inhomogeneous fractional differential equation with continuous variable coefficients. Adv. Differ. Equ., 2019, (2019), \#256.

\bibitem{AMS-1938} E. Pitcher, W.E. Sewel. Existence theorems for solutions of differential equations of non-integer order. Bull. Amer. Math. Soc., 44(2), (1938), 100--107.


\bibitem{RRS} J.E. Restrepo, M. Ruzhansky, D. Suragan. Explicit solutions for linear variable-coefficient fractional differential equations with respect to functions. Appl. Math. Comput. 403, (2021).

\bibitem{RSade} J.E. Restrepo, D. Suragan. Direct and inverse Cauchy problems for generalized space-time fractional differential equations. Adv. Differential Equations 26(7--8) (2021), 305--339.

\bibitem{dirac-joel} J.E. Restrepo, M. Ruzhasky, D. Suragan. Generalized fractional Dirac type operators. Fract. Calc. Appl. Anal. 26, (2023), 2720--2756.

\bibitem{vcapl} M. Rivero, L. Rodr\'igez--Cierm\'a, J. J. Trujillo. Linear fractional differential equations with variable coefficients, Appl. Math. Lett., 21, (2008), 892--897.


\bibitem{samko} S.G. Samko, A.A. Kilbas, O.I. Marichev. Fractional integrals and derivatives, translated from the 1987 Russian original, Gordon and Breach, Yverdon, (1993).

\bibitem{Samko1} S.G. Samko, R.P. Cardoso. Integral equations of the first kind of Sonine type. Int. J. Math. Math. Sci., 2003(57), (2003), 3609--3632.


\end{thebibliography}
\end{document}